\documentclass{article}

\usepackage{arxiv}
\usepackage{hyperref}
\usepackage{lipsum}
\usepackage{amsmath}
\usepackage{amsfonts}
\usepackage{amssymb}
\usepackage{graphicx}
\usepackage{subfigure}
\usepackage{multirow,bigdelim}
\ifpdf
  \DeclareGraphicsExtensions{.eps,.pdf,.png,.jpg}
\else
  \DeclareGraphicsExtensions{.eps}
\fi

\newcommand{\ignore}[1]{}
\usepackage{algorithm}
\usepackage{algpseudocode}

\newtheorem{theorem}{Theorem}

\newcommand{\TheTitle}{High-Performance Partial Spectrum Computation for Symmetric eigenvalue problems and the SVD} 
\newcommand{\TheAuthors}{D. Keyes, H. Ltaief, Y. Nakatsukasa, and D. Sukkari}

\title{\TheTitle
}

\author{
  David Keyes\thanks{King Abdullah University of Science and Technology, Extreme Computing 
  Research Center, Computer, Electrical, and Mathematical Sciences and Engineering Division, 
  Thuwal, 23955-6900, Saudi Arabia
    (\texttt{David.Keyes,Hatem.Ltaief@kaust.edu.sa}).}
	\and
	Hatem Ltaief\footnotemark[1]
	\and
	Yuji Nakatsukasa\thanks{Mathematical Institute, University of Oxford, Oxford OX2 6GG, UK
    (\texttt{Yuji.Nakatsukasa@maths.ox.ac.uk}).}
	\and
	Dalal Sukkari\thanks{Innovative Computing Laboratory, University of Tennessee, Knoxville TN 37996, USA
    (\texttt{sukkari@icl.utk.edu}).}
}

\usepackage{amsopn}
\DeclareMathOperator{\diag}{diag}

\ifpdf
\hypersetup{
  pdftitle={\TheTitle},
  pdfauthor={\TheAuthors}
}
\fi

\begin{document}

\maketitle

\begin{abstract}
Current dense symmetric eigenvalue (EIG) and singular value decomposition (SVD) 
implementations may suffer from the lack of concurrency during the reduction 
step toward the corresponding condensed matrix forms, i.e., tridiagonal and
bidiagonal, respectively. This performance bottleneck is typical for the two-sided
transformations due to the Level-2 BLAS calls. These memory-bound functions
are inherently limited by the speed of the bus bandwidth and may already saturate the memory
bandwidth with only a small number of processes. Therefore,
the current state-of-the-art EIG and SVD implementations
may achieve only a small fraction of the system's sustained peak performance. 
The QR-based Dynamically Weighted Halley (QDWH) algorithm may be used as a pre-processing 
step toward the EIG and SVD solvers, while mitigating the aforementioned bottleneck. 
\texttt{QDWH-EIG} and \texttt{QDWH-SVD} expose more parallelism, 
while relying on compute-bound matrix operations. Both run closer to
the sustained peak performance of the system, but at the expense of performing more 
floating-point operations than the standard EIG and SVD algorithms.
These algorithms are designed to compute the full SVD and eigendecomposition.
In this paper, we introduce a new QDWH-based
solver for computing the partial spectrum for EIG (\texttt{QDWHpartial-EIG}) 
and SVD (\texttt{QDWHpartial-SVD}) problems.
By optimizing the rational function underlying the algorithms only in the desired part of the spectrum,
\texttt{QDWHpartial-EIG} and \texttt{QDWHpartial-SVD} algorithms
efficiently compute a fraction (say $1-20\%$) of the eigenspectrum as well as the most significant singular values/vectors, respectively.
We develop high-performance implementations of \texttt{QDWHpartial-EIG} and \texttt{QDWHpartial-SVD} 
on distributed-memory manycore systems and demonstrate their numerical robustness.  
We perform a benchmarking campaign against their counterparts from the 
state-of-the-art numerical libraries (i.e., \texttt{ScaLAPACK}, \texttt{ELPA}, \texttt{KSVD})
across various matrix sizes using up to $36$K MPI processes. Experimental
results show performance speedups for \texttt{QDWHpartial-SVD} up to $6$X and $2$X 
against \texttt{PDGESVD} from \texttt{ScaLAPACK} and \texttt{KSVD}, respectively.
\texttt{QDWHpartial-EIG} outperforms \texttt{PDSYEVD} from \texttt{ScaLAPACK} up to $3.5$X
but remains slower compared to \texttt{ELPA}. \texttt{QDWHpartial-EIG}
achieves, however, a better occupancy of the underlying hardware by extracting 
higher sustained peak performance than \texttt{ELPA}, which is critical moving forward with 
accelerator-based supercomputers.

\end{abstract}



\section{Introduction}
\label{sec:intro}
Solving the dense symmetric eigenvalue (EIG) and singular value decomposition (SVD)  
problems~\cite{Golub1970,golubbook4th,trefethen1997} represent one of the main computational
phases for many scientific problems, e.g., signal processing~
\cite{soummer2012detection},
pattern recognition~\cite{Elden2007}, statistics~\cite{Oseledets2009}, quantum chemistry~\cite{eig-app1}, 
quantum physics~\cite{eig-app2}, and quantum mechanics~\cite{eig-app3}.
There are actually many applications for which there is interest in 
extracting only a partial eigenspectrum from EIG, e.g., in density function
theory for electronic structure calculations~\cite{Saad1996}. Similarly,
there are several numerical algorithms that may only require the most significant 
singular values with their associated singular vectors from SVD, e.g., determining the 
pseudo-inverse of a matrix~\cite{Ltaief_pasc2} or performing low-rank matrix 
approximations~\cite{akbudak2018hicma,blr-schur,hackbusch2015hierarchical}.

The current state-of-the-art numerical libraries \texttt{LAPACK}~\cite{LAPACK} 
and \texttt{ScaLAPACK}~\cite{ScaLAPACK} for shared-memory
and distributed-memory systems, respectively, provide EIG and SVD implementations.
They first reduce the original dense matrix into condensed tridiagonal and bidiagonal 
forms, before computing the eigenspectrum and the singular values/vectors, respectively. 
Although this initial reduction phase occupies a small part of the floating-point
operations (flops), it may still account for up to half of the overall time taken
by the EIG and SVD solvers. This is due to its memory-bound execution during 
the expensive panel factorization based on Level-2 BLAS, which requires accessing the entire 
unreduced trailing submatrix. The memory bandwidth may quickly become saturated
and adding more computational resources may actually slow down further the 
execution. Two-stage matrix reductions~\cite{Lang99,sbr} for EIG and SVD have become popular,
as they allow to cast some of the Level-2 BLAS operations into compute-bound
Level-3 BLAS.

In many applications, one is interested only in a subset of the spectrum; usually the
extremal (largest/smallest) eigenpairs and the dominant singular triplets. However, the traditional one/two-stage
reduction-based approaches still require to transform the whole matrix
into tridiagonal or bidiagonal form, so the overall runtime is comparable to that of a full decomposition.

In this paper, we design and implement algorithms that remove these
aforementioned limitations in order to compute the partial spectrum 
for the EIG and SVD solvers. Based on the polar decomposition,
these new high performance EIG and SVD algorithms rely on the QR-based Dynamically 
Weighted Halley (QDWH) method to compute a partial spectrum. As initially 
introduced in~\cite{Yuji2010}, QDWH is an expensive approach with
a much higher number of flops when used toward computing the full 
spectrum for EIG and SVD solvers~\cite{Yuji2013}. But ultimately, 
it turns out to be a competitive approach for SVD,
while remaining an interesting alternative for EIG~\cite{Sukkari2016,Sukkari2016bis,Sukkari2017bis}.
The main idea consists in compensating for these extra flops with the higher level of concurrency
and the compute-bound nature achieved by the QDWH numerical kernels.
We then leverage existing QDWH-based SVD and EIG algorithms to compute the 
partial spectrum for the EIG (\texttt{QDWHpartial-EIG}) and SVD (\texttt{QDWHpartial-SVD}).
Our new \texttt{QDWHpartial-EIG} and \texttt{QDWHpartial-SVD} algorithms permit to redirect 
the computational power toward only the operations necessary for 
the computation of the spectrum of interest. This inherent flexibility
of QDWH makes it even more competitive against the existing implementations
for extracting the partial spectrum for the EIG.

We deploy both \texttt{QDWHpartial-EIG} and \texttt{QDWHpartial-SVD} implementations 
on a large distributed-memory system and ensure
the original numerical robustness of QDWH-based full EIG and SVD is maintained, 
for various matrix types. We then assess their individual performance and 
compare them against their respective counterparts from the state-of-the-art 
numerical libraries (i.e., \texttt{ScaLAPACK}~\cite{ScaLAPACK}, \texttt{ELPA}~\cite{elpa}, 
and \texttt{KSVD}~\cite{Sukkari2017bis}). Experimental
results show performance speedups for \texttt{QDWHpartial-SVD} up to $6$X and $2$X 
against \texttt{PDGESVD} from \texttt{ScaLAPACK} and \texttt{KSVD}, respectively.
\texttt{QDWHpartial-EIG} outperforms \texttt{PDSYEVD} from \texttt{ScaLAPACK} up to $2.5$X
while being within reach compared to \texttt{ELPA}. Moreover, \texttt{QDWHpartial-EIG} is capable 
to extract a higher sustained peak performance from the underlying hardware. This is critical 
when looking at hardware architecture trends, where accelerator-based 
systems generously provisioned with flops will constitute most of recently announced exascale supercomputers.

The remainder of the paper is organized as follows. Section~\ref{sec:rw}
describes related work for state-of-the-art EIG and SVD solvers. Section~\ref{sec:bg}
reviews the background of the QDWH approach for the polar decomposition
and its application to EIG and SVD solvers. Section~\ref{sec:partial}
introduces the new \texttt{QDWHpartial-EIG} and \texttt{QDWHpartial-SVD} algorithms for
computing only the partial spectrum.
Section~\ref{sec:impl} describes the implementation details and
Section~\ref{sec:complexity} estimates the algorithmic operation counts.
Section~\ref{sec:acc} highlights the numerical robustness of
both \texttt{QDWHpartial-EIG} and \texttt{QDWHpartial-SVD} algorithms.
Section~\ref{sec:results} assesses the achieved performance results and we conclude in
Section~\ref{sec:conclusion}.
\section{Related work}
\label{sec:rw}
When computing the \emph{full} spectrum for symmetric (or Hermitian for complex matrices) 
eigenvalue and SVD solvers, the state-of-the-art approaches can be split 
into two categories. The one-stage approaches, as implemented in \texttt{LAPACK}~\cite{LAPACK} 
and \texttt{ScaLAPACK}~\cite{ScaLAPACK}, reduce the dense matrix into a condensed form using
a single phase of orthogonal transformations, before extracting the spectrum of interest.
To further promote Level-3 BLAS operations during this single stage, two-stage approaches~\cite{Lang99,sbr} 
have emerged as an efficient algorithmic alternative in better extracting the 
hardware performance. Although they come at the price of extra floating-point
operations (flops), their high performance implementations have contributed in their
wide adaption in the software ecosystem 
within the \texttt{PLASMA}~\cite{luszczek2011ipdps,Haidar11,haidar1,Ltaiefparco11,ltaief-tilecasvd,ltaief2} 
and \texttt{MAGMA}~\cite{HaidarTDSS14} libraries on shared-memory systems (possibly equipped with GPUs) for
eigenvalue and SVD solvers
or the \texttt{ELPA}~\cite{elpa} and \texttt{EigenExa}~\cite{eigenexa} libraries for 
only eigenvalue solvers on distributed-memory systems. 

When it comes to calculating the \emph{partial} spectrum for eigenvalue and SVD solvers,
the one and two-stage approaches are inefficient as described above. 
A more recent work~\cite{Marques2020}
shows how to partially compute the SVD out of the bidiagonal form using an associated tridiagonal eigenproblem.
But yet again, the condensed form remains the ultimate starting point and one of the most expensive computational operations.

Completely different classes of algorithms for computing a small part of the spectrum have been developed,
most prominently the Lanczos algorithm (more generally Krylov subspace methods)  and the more recent
randomized algorithms~\cite{halko2011finding}. While these can be very  powerful, they come with certain
drawbacks in the situation that we consider. 

Krylov methods are usually suitable only when a very small fraction of the spectrum (usually $O(1)$
eigenvalues or singular values) is required, and sometimes fails to provide full accuracy. In this
work we consider the case where a nonnegligible portion of the spectrum (say $1-20\%$) is desired. 

Randomized algorithms can be an extremely effective means of finding an approximate SVD, and are
rapidly gaining popularity. However, they usually come with poorer accuracy guarantees, giving outputs that are suboptimal by an $O(1)$ factor;
see~\cite[\S 10]{halko2011finding}, \cite[\S 3]{nakatsukasa2020fast} (these guarantees are
still remarkable---especially when the spectrum decays rapidly---and enough in many applications~\cite{akbudak2018hicma}). 

In this paper, we propose to revisit and modify the \texttt{QDWH-EIG/SVD} algorithms~\cite{Yuji2010,Yuji2013}
in order to provide support for determining only a partial spectrum for the eigenvalue and
SVD solvers. We aim to compute the eigen/singular values and vectors essentially to full working precision. 
These algorithms and their high performance
implementations~\cite{Sukkari2016bis,Sukkari2016,Sukkari2017bis,Sukkaritpds} do not require
a reduction to tridiagonal or bidiagonal forms. They iteratively compute the
polar decomposition---based on conventional, compute-bound, and highly-parallel dense linear algebra 
operations---as a preprocessing step toward the eigenvalue and SVD solvers. By altering the core
algorithmic feature of \texttt{QDWH-EIG/SVD}, the new \texttt{QDWHpartial-EIG/SVD} approach transforms directly the original 
dense matrix to a much smaller one, with a size roughly of the spectrum of interest. Since the transformation
occurs at the beginning of the \texttt{QDWHpartial-EIG/SVD} procedure,  
the power of computational resources is tailored solely to operations that are intimately related to the 
eigenspace of interest.
\section{QDWH-based polar decomposition and its application
to full symmetric eigenvalue and SVD solvers}
\label{sec:bg}
The Polar Decomposition (PD) $A=U_pH\in\mathbb{C}^{m\times n}$, where
$U_p\in\mathbb{C}^{m\times n}$ is the \emph{unitary polar factor}
with $U_p^*U_p=I_n$ and $H$ is Hermitian positive semidefinite, exists for any matrix.
It is an important
matrix decomposition for various applications, including inertial navigation~\cite{polar-app1},
chemistry~\cite{polar-app2}, and computation of block reflectors in numerical linear algebra~\cite{polar-app3}.
It can be used as a first computational phase toward computing the EIG/SVD~\cite{Yuji2013} in the context
of the QR-based Dynamically Weighted Halley (QDWH) method.

\subsection{The QDWH-Based PD Algorithm}
\label{subsec:qdwhpd}
The dynamically weighted Halley iteration to find the PD can be
summarized as follows:
\begin{equation}
\begin{aligned}
X_{0} &= A/\alpha,\\
X_{k+1} &= X_{k} (a_{k}I + b_{k} X_{k}^* X_{k})(I + c_{k}X_{k}^* X_{k})^{-1}.
\label{polar-itr}
\end{aligned}
\end{equation}
The scalars ($a_{k}, b_{k}, c_{k}$) are chosen dynamically to speed up 
the convergence~\cite{Yuji2010}.
More specifically, they are chosen so that the rational function
$r_k(x)=x(a_k+b_kx^2)/(1+c_kx^2)$ is the scaled \emph{Zolotarev} function of type $(3,2)$,
the best rational approximation to the sign function on $[-1,-\ell_k]\cup [\ell_k,1]$.
Here $\ell_0=1/\kappa_2(A)$ (or its estimate) and follows the updating formula $\ell_{k}=r_{k}(\ell_{k-1})$.
The singular values of $X_{k}$ are $\Sigma_{k} = r_{k}( \cdots r_{2}(r_{1}(\Sigma)))$, and lie in $[\ell_k,1]$.
Remarkably, the composition of the rational functions $r_{k}( \cdots r_{2}(r_{1}(\Sigma))))$ is 
again a Zolotarev function, of much higher type $(3^k,3^k-1)$. 
Together with the exponential convergence of Zolotarev functions, 
QDWH converges in at most \emph{six} iterations to obtain $X_k\rightarrow U_p$
(and $\ell_k\rightarrow 1$) in double precision for matrices with $\kappa_2(A)\leq 10^{15}$. 

Based on the fact~\cite[p.~219]{Higham08} that 
$c X (I + c^{2} X^*X)^{-1}=Q_{1}Q_{2}^* $, where 
$\begin{bmatrix}
c X \\
I
\end{bmatrix}
= \begin{bmatrix}
Q_{1} \\
Q_{2}
\end{bmatrix}R
$
is the $QR$ decomposition,  with $X, Q \in \mathbb{R}^{m \times n}$ and $Q_{2}, R \in \mathbb{R}^{n \times n}$, 
the Equation~\eqref{polar-itr} can be replaced with the following 
inverse-free and stable $QR$-based implementation~\cite{Yuji2013}:
\begin{equation}
\begin{aligned}
\begin{bmatrix} \sqrt{c_{k}}X_{k} \\ I  \end{bmatrix} &= 
\begin{bmatrix} Q_{1} \\ Q_{2}  \end{bmatrix}R,\\ 
X_{k+1} &= \frac{b_{k}}{c_{k}}X_{k} + \frac{1}{\sqrt{c_{k}}} \left( a_{k} - \frac{b_{k}}{c_{k}} \right) Q_{1}Q_{2}^*. 
\label{practical-polar-qr}
\end{aligned}
\end{equation}
This 
 represents the QR-based Dynamically                                                  
Weighted Halley (QDWH) algorithm. Further details can be found in~\cite{Yuji2016}.

After a few QDWH iterations from Equation~\eqref{practical-polar-qr}, the $X_{k}$ eventually becomes well-conditioned
and $\ell_k=O(1)$. Once this happens, a lower-cost Cholesky-based iteration can be used instead, as follows:
\begin{equation}
\begin{aligned}
X_{k+1} &= \frac{b_{k}}{c_{k}}X_{k} +  \left( a_{k} - \frac{b_{k}}{c_{k}} \right) (X_{k}W_{k}^{-1})W_{k}^{-*},\\
W_{k} &= \mbox{chol}(Z_{k}),\:Z_{k} = I + c_{k}X_{k}^* X_{k}.
\label{practical-polar-chol}
\end{aligned}
\end{equation}

A higher-order variant of the QDWH algorithm that employ higher-degree Zolotarev functions has been developed~\cite{Yuji2016},
which further increases the degree of parallelism.


\subsection{Applying QDWH to Full Symmetric Eigenvalue and SVD Solvers}
\label{subsec:eigsvd}


First, we recall the mechanism of \texttt{QDWH-EIG}~\cite{Yuji2013}, on which the new \texttt{QDWHpartial-EIG} 
algorithm will be based. 
Let $A$ be an $n\times n$ symmetric matrix and write 
\begin{align}
  A &= V \diag(\Lambda_+, \Lambda_-) V^* \nonumber\\
    &= V \diag(I_{n-k}, -I_{k}) V^* \cdot V  \diag(\Lambda_+, |\Lambda_-|) V^* \nonumber\\
    &\equiv U_p H,   \label{eq:pdeig}
\end{align}
where $U_pH$ is the polar decomposition~\cite[Ch.~8]{Higham08}.
$k$ is the number of negative eigenvalues, which we do not assume to be known. 
Suppose that the unitary polar factor $U_p$ is computed. 
This means we have mapped all the eigenvalues to $1$ or $-1$. 
We partition $V = [V_+,~V_-]$ conformably with $\Lambda$, and 
note that
\begin{align}
	\frac{1}{2}(I-U_p) &= \frac{1}{2}\left(I-[V_+~V_-]
\begin{bmatrix}
I_{n-1}&0\\ 0&-I_{k}  
\end{bmatrix}
[V_+\ V_-]^*  \right)\nonumber\\
        &= [V_+~V_-]
        \begin{bmatrix}
0&0\\ 0&I_{k}
        \end{bmatrix}[V_+\ V_-]^*\nonumber\\
		 &= V_-V_-^*.
\end{align}
Hence the
symmetric matrix $C=\frac{1}{2}(I-U_p)=V_-V_-^* $ is an orthogonal
projector onto $\mbox{Span}(V_-)$, the invariant subspace corresponding to the 
negative eigenvalues. 
We can then project the matrix $A$ (Rayleigh-Ritz process)
to obtain the eigenvalues and eigenvectors: 
 the eigenvalues of $V_-^*AV_-$ are equal to those of $\Lambda_-$, and 
denoting by $V_-^*AV_-=W\Lambda_-W^*$  the eigenvalue decomposition, 
we see that $V_-W$ is the matrix of eigenvectors.
Analogously, we can obtain $V_+$ by finding the subspace spanned by 
$ \frac{1}{2}(I+U_p)$. 

Second, the polar decomposition can be also used directly
toward calculating the SVD, i.e., 
$A\; =\;U_{p}H\; =\; U_{p} (V\Sigma V^{*}) \;= \; (U_{p} V)\Sigma V^{\top} \; = \;U\Sigma V^{*}$
where $\Sigma$ is 
the matrix containing all the singular values, 
and $U$ and $V$ are the orthogonal matrices containing the left and 
right singular vectors, respectively. The resulting \texttt{QDWH-SVD} procedure
relies on \texttt{QDWH-EIG} (or any other eigensolvers) to compute the intermediate 
eigendecomposition $H=V\Sigma V^*$ required for the final SVD.


\section{Leveraging QDWH for computing the partial spectrum
of symmetric eigenvalue and SVD solvers}
\label{sec:partial}
In this section, we present modified versions of \texttt{QDWH-EIG} and \texttt{QDWH-SVD} to compute
the partial (negative) eigenspectrum of a symmetric matrix (\texttt{QDWHpartial-EIG})
and to extract the most significant singular values 
and their corresponding singular vectors (\texttt{QDWHpartial-SVD}). 

In what follows we treat nonreal matrices $A\in\mathbb{C}^{n\times n}$; when $A$ is real, the superscripts $\cdot^*$ should be replaced by $\cdot^\top$ and everything can be executed using only real arithmetic. 

\subsection{QDWHpartial-EIG}
\label{subsec:ses}
We introduce \texttt{QDWHpartial-EIG}, an algorithm for computing the 
negative eigenvalues and its corresponding eigenvectors 
(of course the algorithm can be modified trivially to find the 
eigenvalues/vectors smaller/larger than any specific number 
by shifting and scaling by $-1$).

For simplicity of exposition, here we suppose that the negative eigenvalues lie (roughly) in $[-1,0)$, that is, 
$\lambda_{\min}(A)\gtrsim -1$. To ensure this we need a lower bound $\mu\leq \lambda_{\min}(A)<0$; many algorithms are available for this task; we use a few steps of the Lanczos iteration to estimate $\lambda_{\min}(A)$. We then scale the matrix $A:=A/|\mu|$. 

Recall that the mathematics underlying the \texttt{QDWH-EIG} algorithm
is rational approximation: it 
finds a rational function $r$ that approximates the sign function, so that it 
maps the negative eigenvalues to $-1$, and positive eigenvalues to $1$. 
Thus $r(A)$ has eigenvalues $\pm 1$.

Now, suppose that we only need the negative eigenspace $V_-$. Can we cut corners? The answer 
is yes---and this is not just that we can skip computing $V_*$ 
once $U_p$ is computed. We shall avoid computing $U_p$. Essentially, 
we need to ensure only that the negative eigenvalues have been 
mapped to $-1$; the positive eigenvalues are irrelevant. 
Namely, the idea is to find a matrix whose null space $Q_2$ contains $V_-$. We then extract $V_-$ as a subspace of $Q_2$. 

To obtain such $Q_2$ efficiently, a key idea is to only map the required eigenvalues by the sign function. 
QDWH works on the interval $[-1,-\ell_0]\cup[\ell_0,1]$ and approximates the sign function there with a rational function $r$; the QDWH convergence is governed by $\ell_0$; the larger the faster. 

Now, we can use the \emph{shifted} rational function $r(x-s)$ for some  $s\in(0,1)$; if $r$ approximates the sign function on $[-1,-s]\cup[s,1]$ (taking $\ell_0=s$), then we have $r(x-s)\approx -1$  on $[-1+s,0]$. 
Therefore, defining 
$\tilde A := (1-s)A-sI$
(whose negative eigenvalues are again in $[-1,0]$, and the negative eigenvalues of $A$ correspond to those of $\tilde A$ in $[-1,-s]$), 
we see that the negative eigenvalues of $A$ are mapped to $-1$ by 
$r(\tilde A)$. 

\begin{figure}[htbp]
\centering
      \includegraphics[width=0.4\textwidth]{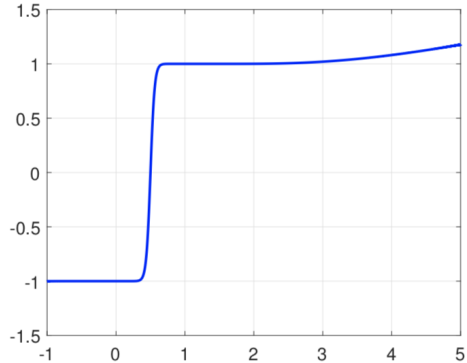}  
  \caption{A type $(9,8)$ rational function $y=r(x-s)$ that 
maps the interval $[-1,0]$ to $-1$. $s=0.5$. 
}
  \label{fig:rx}
\end{figure}

The reason we introduce the shift $s$ is that we can then set $\ell_0=s$, so QDWH can converge faster with a low-degree Zolotarev function. 
For illustration, Fig.~\ref{fig:rx} shows a typical plot of $r(x-s)$. 
Observe  how the interval $[-1,0]$ is mapped to $-1$ with a low-degree rational function.


Suppose that we have computed 
the matrix function $r(\tilde A)$. It has 
$k$ (or more) eigenvalues at $-1$, 
where $k$ is the number of negative eigenvalues of $A$. 
Since the desired eigenvalues have been mapped to $-1$, we are interested in the null space of $r(\tilde A)+I$. 
While $r(\tilde A)$ has $n-k$ (or fewer) eigenvalues away from $-1$, the knowledge of their precise values
is not critical. It does matter, however, that they do not blow up to $\pm\infty$, for numerical stability.

It is worth emphasizing that there is a 
nontrivial interval on the positive side (on $[0,\hat s]$ for $\hat s<s$) that is mapped close to $-1$. 
This is an inevitable nature of rational functions, and can cause the dimension of the numerical null space of the matrix $r(\tilde A)+I$ 
to be larger than $k$. Indeed if we take $s$ too large to a point where all the eigenvalues 
of $A$ are mapped close to $-1$, then 
$r(\tilde A)+I$ converges to zero, and the process below 
leads to no efficiency gain 
(i.e., the whole space becomes the computed null space and no gain is 
obtained relative to doing a full eigendecomposition). 
We discuss how to choose an appropriate $s$ shortly. 

We then work with the matrix $r(\tilde A)+I$. This matrix 
is not a partial isometry 
as it was in QDWH, but it is still rank deficient, with deficiency $k$ or more. 

To compute the null space of $r(\tilde A)+I$ (or a larger subspace that contains it), we compute the QR factorization 
\begin{equation}  \label{eq:qrfac}
r(\tilde A)+I = [Q_1,Q_2]
\begin{bmatrix}
R_{11}  & R_{12}\\
0 & R_{22}
\end{bmatrix}  
\end{equation}
where the size of $R_{11}$ (and hence of $R_{22}\in\mathbb{C}^{\ell\times\ell}$, where $\ell\geq k$) is chosen so that it only has ``large" singular values; the idea is that $Q_2$ then contains the null space of $A$ that we require. 
To quantify the claim we use matrix perturbation theory. 
\begin{theorem}\label{thm:pert}
Let $B = [Q_1,Q_2]\begin{bmatrix}
R_{11}  & R_{12}\\
0 & R_{22}
\end{bmatrix} \in\mathbb{C}^{m\times n} (m\leq n)$ be a full QR factorization such that $[Q_1,Q_2]\in\mathbb{C}^{m\times m}$ is unitary with 
$Q_2\in\mathbb{C}^{m\times \ell}$ and 
$R_{11}\in\mathbb{C}^{(m-\ell)\times (m-\ell)}$. 
Let $V_0\in\mathbb{C}^{m\times k}$ have orthonormal columns $V_0^*V_0=I_k$, with $k\leq \ell$. Then 
\begin{equation}  \label{eq:threshold}
\sin\angle(V_0,Q_2)\leq \frac{\|V_0^*B\|_2}{\sigma_{\min}(R_{11})}.
\end{equation}
\end{theorem}
{\it Remark}. 
The main situation of interest is when $V_0$ spans 
an approximate left null space of $B$ such that $\|V_0^*B\|_2=O(u\|B\|_2)$, where $u$ is unit roundoff; then the theorem shows the subspace $V_0$ is approximately contained in $Q_2$, up to $O(u/\sigma_{\min}(R_{11}))$. 

{\sc proof}. 
First recall that 
the canonical angles 
$\angle_1(V_0,Q_1),\ldots,\angle_k(V_0,Q_1)$
between two subspaces 
of dimensions $k,\ell(\ell\geq k)$
spanned by the orthonormal matrices 
$V_0\in\mathbb{C}^{m\times k}$ and $Q_2\in\mathbb{C}^{m\times \ell}$ (for which $Q_1$ is the orthogonal complement $Q_1^*Q_2=0$) are defined by
$\sin\angle_i(V_0,Q_1)=\sigma_i(V_0^*Q_1)$~\cite[\S~2]{golubbook4th} for $i=1,\ldots,k$. 
It thus suffices to show that $\|V_0^*Q_1\|_2\leq \frac{\|V_0^*B\|_2}{\sigma_{\min}(R_{11})}$. 

Now we have 
\[V_0^*B=[V_0^*Q_1,V_0^*Q_2]\begin{bmatrix}
R_{11}  & R_{12}\\0 & R_{22}\end{bmatrix}  =\begin{bmatrix}V_0^*Q_1R_{11},V_0^*Q_1R_{12}+V_0^*Q_2R_{22}\end{bmatrix}.\]
Hence we obtain $\|V_0^*B\|_2\geq \|V_0^*Q_1R_{11}\|_2\geq \|V_0^*Q_1\|_2\sigma_{\min}(R_{11})$. 
It follows that 
$\|V_0^*Q_1\|_2\leq \frac{\|V_0^*B\|_2}{\sigma_{\min}(R_{11})}$, as required. 
\hfill$\square$

Let us make two remarks about the theorem. 
\begin{itemize}
\item Theorem~\ref{thm:pert} does not require $B$ to be symmetric, or even square; $m<n$ is allowed. It does not apply directly to the $m>n$ case, as the left null space $V_0$ is larger than the rank deficiency of $B$. 
\item Finding a numerical null space of a matrix is a classical problem in numerical linear algebra, and a reliable algorithm is usually based on either the SVD or a strong rank-revealing QR factorization~\cite{gu1996efficient}. The assumptions in the theorem are much weaker; the reason they suffice is that~\eqref{eq:threshold} only states that the null space $V_0$ is \emph{contained} in (and not equal to) $Q_2$; in particular, it does not claim $\|BQ_2\|_2$ is small. Once such $Q_2$ is obtained, our algorithm will extract a null space of $BQ_2$ using the Rayleigh-Ritz process. 
\end{itemize}
When Theorem~\ref{thm:pert} is applied with
$B\leftarrow r(\tilde A)+I$ as in~\eqref{eq:qrfac} (so $m=n$), it gives
$\sin\angle(V_0,Q_2)\leq \frac{\|V_0^*(r(\tilde A)+I)\|_2}{\sigma_{\min}(R_{11})}$, where 
$V_0$ is the numerical null space of $r(\tilde A)+I$, which has dimension $\geq k$ by construction. 
Therefore $\|V_0^*(r(\tilde A)+I)\|_2=O(u)$, and it follows that 
we have $\sin\angle(V_0,Q_2)\leq O(u/\sigma_{\min}(R_{11}))$, so $Q_2$ contains $V_0$ to working accuracy provided that $\sigma_{\min}(R_{11})\geq tol $ for some tolerance $tol=\Omega(1)$, for example $tol= 0.01$, which is the choice we make by default. Choosing it too large, $tol\approx 1$, results in $\ell\approx n$ so no computational savings, while $tol\ll 1$ causes loss of accuracy. 

\paragraph{Choosing the subspace size $\ell$}
We have seen that by looking for an $\ell$ such that $R_{11}\in\mathbb{C}^{(n-\ell)\times (n-\ell)}$ satisfies the condition 
$\sigma_{\min}(R_{11}) \gtrsim tol$, 
we can find a subspace $Q_2\in\mathbb{C}^{n\times \ell}$ that contains the desired subspace. We would like to find the smallest possible $\ell$ to reduce the cost of the subsequent operations. 

Fortunately, the condition $\sigma_{\min}(R_{11}) \gtrsim tol$ can be checked reliably without computing the singular values of $R_{11}$. 
The key fact is that if~\eqref{eq:qrfac} is a rank-revealing QR factorization, then the $i$th diagonal element of $R_{11}$ is a good approximation to $\sigma_i(R_{11})$~\cite{martinsson2015blocked}. This more than suffices, given that violation by an $O(1)$ factor in the condition $\sigma_{\min}(R_{11}) \gtrsim tol$ only reduces the final accuracy by an $O(1)$ factor. 
It follows that we can simply examine the diagonal entries of $R$ in the QR factorization of $r(\tilde A)+I$ to determine $\ell$. 

In many cases, the QR factorization without pivoting is already rank revealing; if not, or to ensure this is true with high probability, one can take
the QR of 
$(r(\tilde A)+I)\Omega$ for an $n\times n$ random Gaussian matrix $\Omega$--it is known that such QR factorization is rank revealing with high
probability~\cite{Ballard}. (We can alternatively use QR with pivots, but this incurs substantial communication overhead in parallel computing
settings.) 

Once we have computed such a QR factorization, we 
look for the first diagonal $(i,i)$-element of $R$ that comes below $tol$, and take $i-1$  to be the size of $R_{11}$, 
that is, $\ell = n-i+1$. 


Once $Q_2$ is obtained, we then extract the desired space $V_0$ from it by the Rayleigh-Ritz process: compute the the eigenvalues of the $\ell\times\ell$ matrix $Q_2^*AQ_2$, whose negative eigenvalues should match 
those of $\Lambda_-$. 
Denoting by $Q_2^*AQ_2=W\Lambda_-W^*$  the eigenvalue decomposition, 
the $Q_2W$ is the matrix of eigenvectors. 
By construction, $Q_2^*AQ_2$ usually contains eigenvalues that 
are positive, and we discard those.
 $W$ is then an $\ell\times k$ 
matrix corresponding to the negative eigenvalues. 

We note that since we are interested in 
extremal eigenvalues, Rayleigh-Ritz is a reliable means to extract the 
desired subspace.


\paragraph{Choosing $s$}
Let us consider in more detail the choice of the shift parameter $s$. 
The qualitative behavior has been explained already; taking $s$ small 
results in the Zolotarev function being a poor approximant 
to $\mbox{sign}(x)$ in the interval $[-1,0]$ that we care about, 
while a large $s$ results in good approximation 
$r(x)\approx -1$ on $[-1,0]$, but (undesirably) also on a significant 
positive interval, resulting in the projected size being large. 

\begin{figure}[htpb]
  \begin{minipage}[t]{0.49\hsize}
      \includegraphics[width=0.95\textwidth]{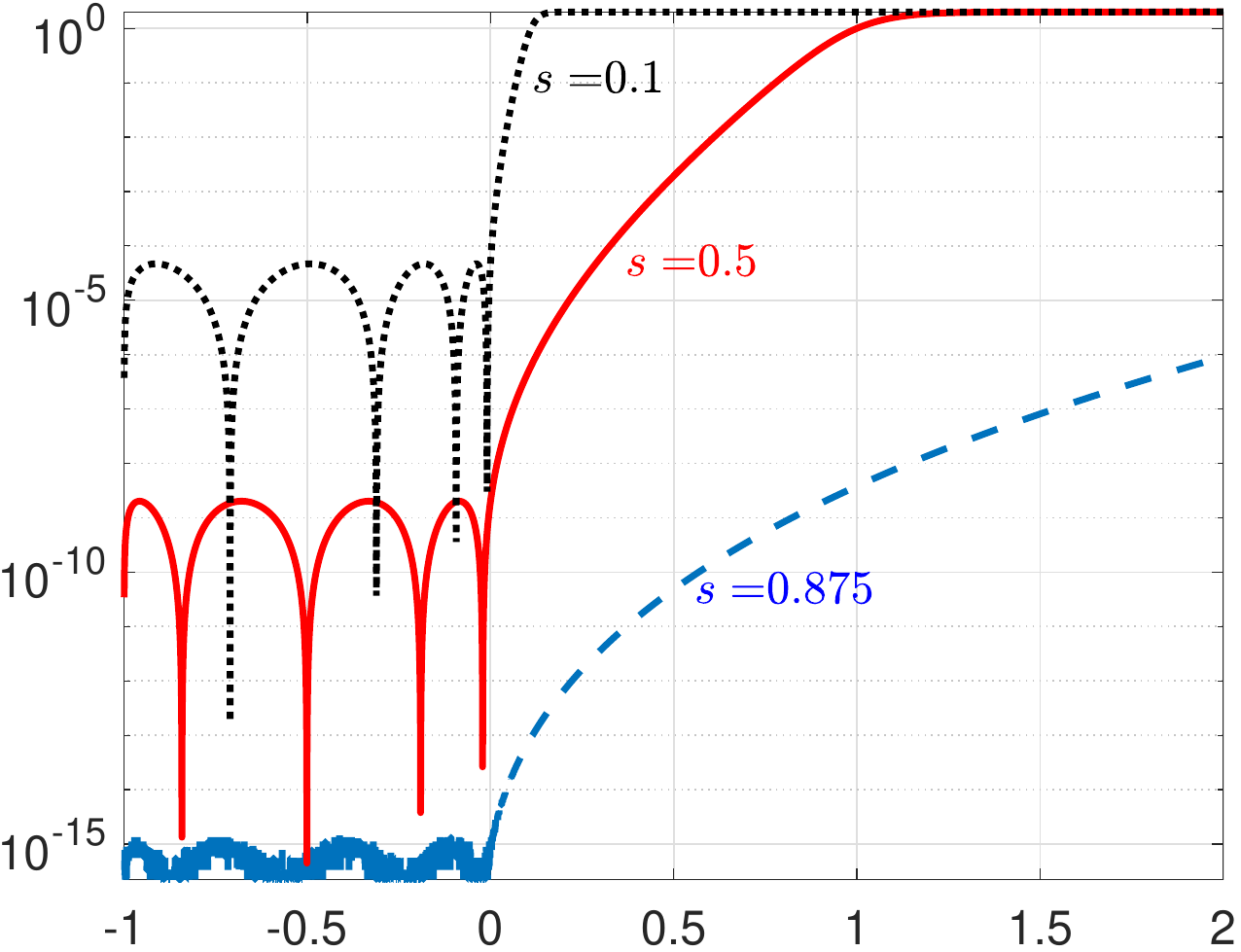}  
  \end{minipage}   
  \begin{minipage}[t]{0.49\hsize}
      \includegraphics[width=0.95\textwidth]{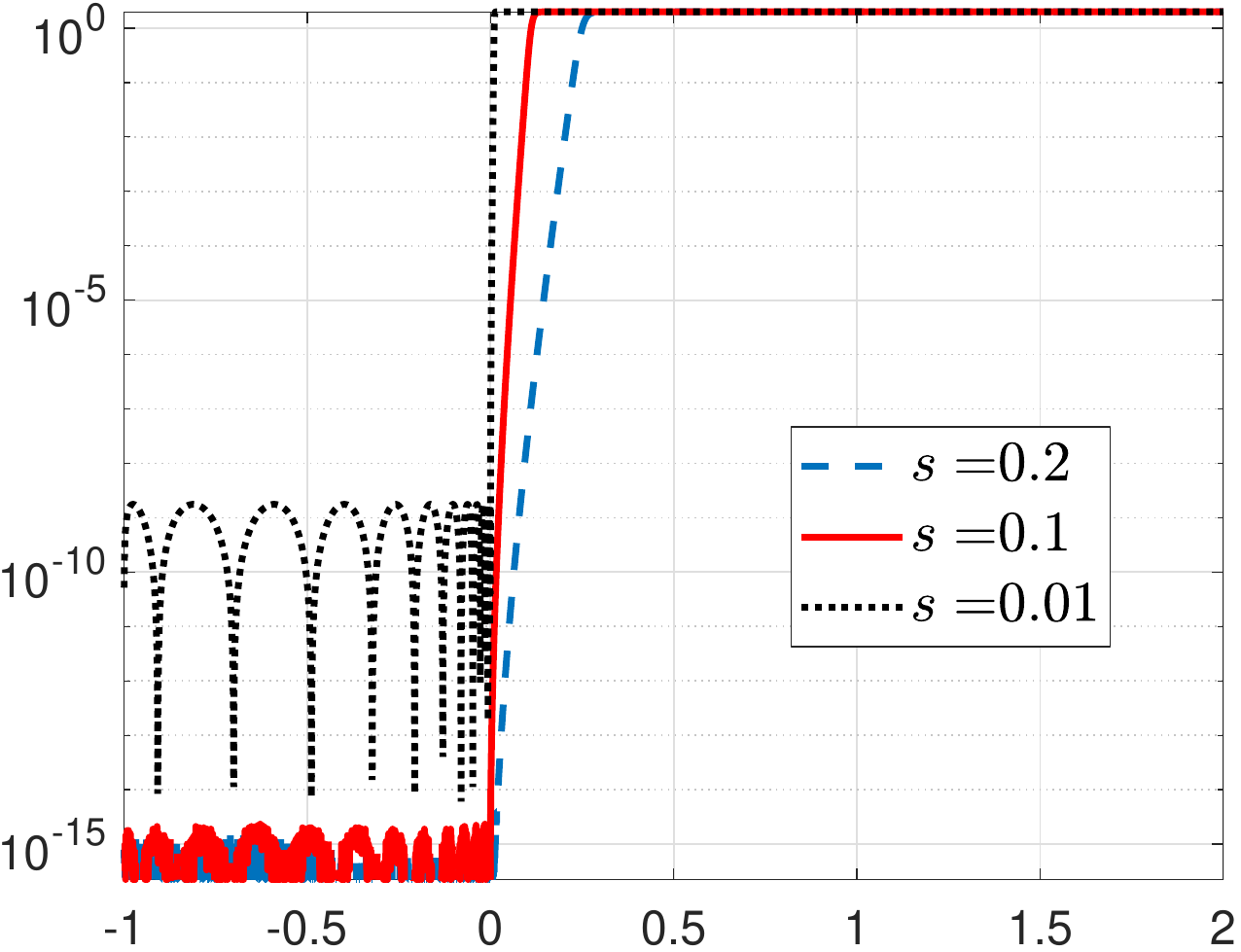}  
  \end{minipage}
  \caption{Semilogy plot of $r(x-s)+1$ for different choices of $s$, 
  for Zolotarev functions of type $(9,8)$ (left) and $(27,26)$ (right). 
}
  \label{fig:ratmm}
\end{figure}

We illustrate this in Fig.~\ref{fig:ratmm}. 
For example, the left plot uses two QDWH iterations, so a Zolotarev function of type $(9,8)$. 
Taking $s=0.875$ gives $O(u)$ approximation 
quality on $[-1,0]$, but the region in which $r(x)+1$ 
is close to $0$ extends far into the positive axis. Taking $s=0.1$, on the 
other hand, solves that issue, but the approximation 
quality on $[-1,0]$ is evidently worse. A similar behavior is seen on the 
right plot, where we use three QDWH iterations and hence a Zolotarev function of higher type $(27,26)$. 
Here the essence stays the same, but a much smaller $s$ is enough to obtain 
$O(u)$ accuracy on $[-1,0]$, and the functions $r(x)$ grow 
steeply until $\approx 1$ for $x>0$. 

Since the goal of the \texttt{QDWHpartial-EIG} iterations is to map the eigenvalues in 
$[-1,0]$ to $0$ (to working precision $O(u)$) by the rational 
function $r(x)+1$ while keeping the positive eigenvalues well separated 
from $0$, the above observation leads to the following strategy 
for choosing $s$:
\begin{enumerate}
\item Determine the type $(2m+1,2m)$ of rational function $r(x)$ to be used. 
\item Choose the smallest $s$ so that $|r(x)+1|\leq O(u)$ on $[-1,0]$. 
\end{enumerate}
In practice, the type $(2m+1,2m)$ is chosen depending on the computational 
budget, and we shall mainly focus on two values $2m+1=3^2$
and $2m+1=3^3$, as these lend to particularly efficient evaluation, 
by taking advantage of the optimality of Zolotarev functions 
under composition~\cite{Yuji2016}. Specifically, they correspond to 
taking two $(2m+1=3^2)$ and three $(2m+1=3^3)$ QDWH iterations. 
For each choice, the value of $s$ satisfying the second condition 
above is found (by simple experiments) to be 
\begin{enumerate}
\item Type $(9,8)$, two QDWH iterations: $s=0.875$,
\item Type $(27,26)$, three QDWH iterations: $s=0.2$. 
\end{enumerate}

These choices are shown in the two plots in Fig.~\ref{fig:ratmm}. 
In most cases, taking three QDWH iterations is recommended as the overhead 
is not too much, while the benefit is significant, 
as can be seen clearly in Fig.~\ref{fig:ratmm}: the function $r(x)+1$ 
of the corresponding cases ($s=.875$ and $s=0.2$) 
take $O(1)$ values for $x\gtrsim 0.1$ with three QDWH iterations, 
but with two iterations, we require $x\gtrsim 5$. 
Consequently,  the projected matrix size $\ell$ will be approximately 
equal to the number of eigenvalues in $[-1,0.1]$ with three 
QDWH, but with two iterations, it captures eigenvalues in $[-1,5]$, 
rendering \texttt{QDWHpartial-EIG} useless unless there is a significant 
portion of large and positive eigenvalues in $A$. We therefore choose the default to be $s=0.2$ and three QDWH iterations. Note that in all iterations, one can safely use the fast Cholesky-based implementation in Eq.~\eqref{practical-polar-chol}. 
Algorithm~\ref{algo:partial-eig} presents the main computational steps 
of \texttt{QDWHpartial-EIG}.

\begin{algorithm}
	\footnotesize{
\begin{algorithmic}[1] 
\State Find an approximate lower bound $\mu\lesssim \lambda_{\min}(A)(<0)$ using a few steps of the Lanczos algorithm, and set $A:=A/|\mu|$. 
\State Set $s=0.2$ (or $s=0.875$) and apply three (or two) iterations of QDWH (using Eq.~\eqref{practical-polar-chol}) to $\tilde A:=(1-s)A-sI$ with $\ell_0=s$, to obtain $r(\tilde A)$ where $r$ is of type $(3^3,3^3-1)$ (or type $(3^2,3^2-1)$). 
\State Compute the QR factorization $\frac{1}{2}(r(\tilde A)+I)=QR$ 
(or $\frac{1}{2}(r(\tilde A)+I)\Omega=QR$, where $\Omega$ is a Gaussian matrix). 
\State Find the index $ind = \min(find(abs(diag(R))<tol=0.01))$.
\State Extract the final columns of $Q$ to get $Q_2 = Q(:,ind:end)$.
\State Rayleigh-Ritz: Compute $Q_2^* A Q_2$ and its eigenvalue decomposition $Q_2^* A Q_2 = \tilde V \Lambda \tilde{V}^*$. 
Extract the parts corresponding to the negative eigenvalues, 
$\Lambda_-,\tilde V_-$. 
\State Rescale $\Lambda:=|\mu|\Lambda_-$ and 
output $\Lambda$ (negative eigenvalues) and $V: = Q_2 \tilde V_-$ (eigenvectors) such that $AV\approx V\Lambda$. 
\end{algorithmic}
}
\caption{\texttt{QDWHpartial-EIG}. Given a symmetric or Hermitian matrix $A$, computes the negative eigenvalues and corresponding eigenvectors. 
}
\label{algo:partial-eig}
\end{algorithm}

Clearly, the algorithm is able to compute the positive eigenvalues, or those that are smaller or larger than a prescribed value by working with a shifted matrix $A-cI$.

\subsection{QDWHpartial-SVD}
\label{subsec:svd}
Given a general matrix $A\in\mathbb{C}^{m\times n} (m\geq n)$, we next consider the task of computing its dominant singular triplets, namely 
computing the singular values and singular vectors corresponding to the 
singular values above a given user-specified relative threshold $s>0$; we use the same letter as the shift for \texttt{QDWHpartial-EIG} as they play a similar role.

As we shall see, essentially the same idea can be applied of optimizing the rational function only in the desired part. 
It is nonetheless worth noting that here we cannot use shifts, as 
to shift singular values without modifying the singular vectors we need the unitary polar factor, which is expensive to compute. 

\begin{figure}[htbp]
\centering
      \includegraphics[width=0.4\textwidth]{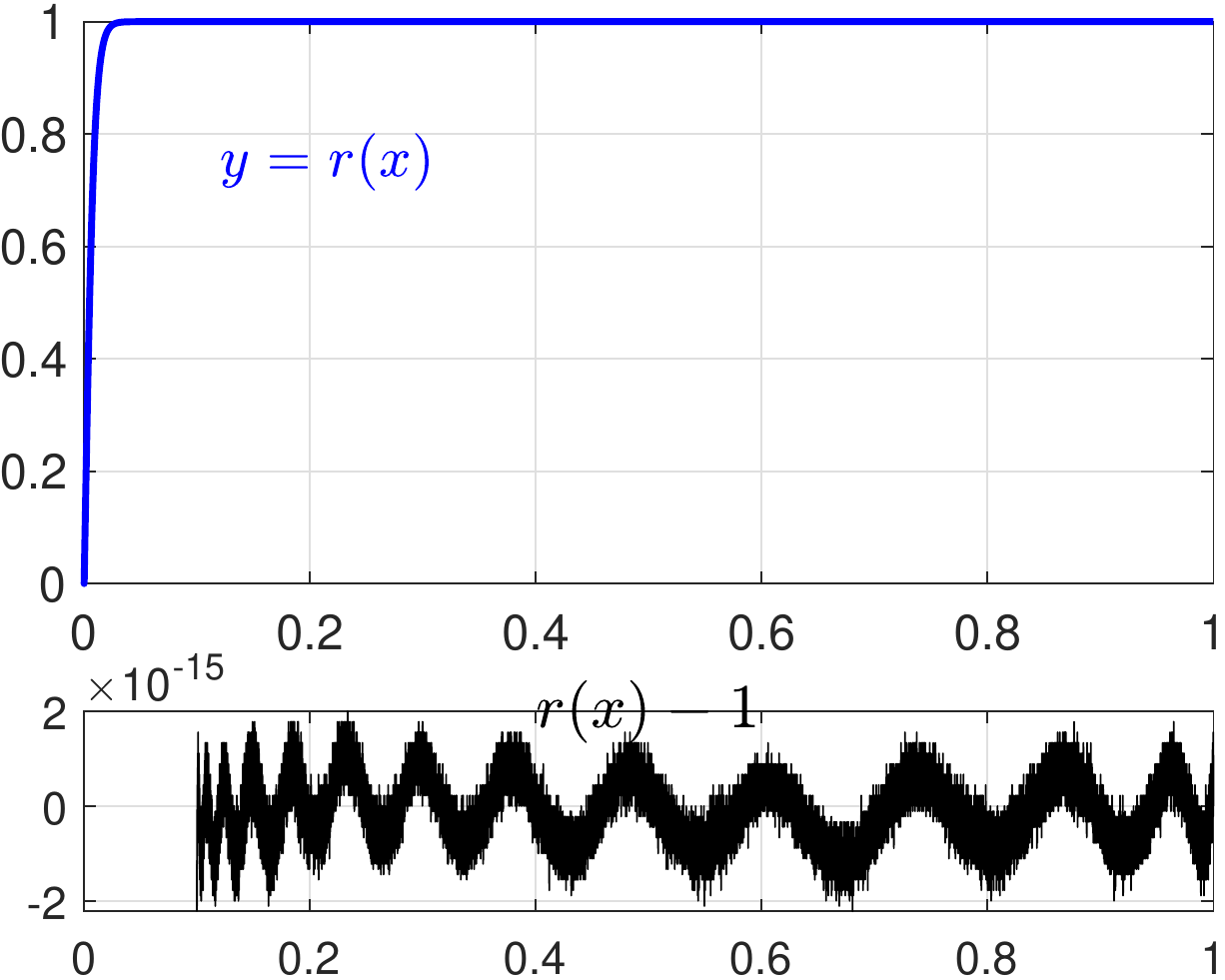}  
  \caption{A type $(27,26)$ Zolotarev function $r$ (top) that 
maps the interval $[s,1]$ to $1$ for $s=0.1$, and its the error $r-1$ (below).
Note the error is $O(u)$ on $[s,1]$; illustrating that three QDWH iterations is enough to find the singular values larger than $s=0.1$. 
}
  \label{fig:rxsvd}
\end{figure}

Assuming w.l.o.g. that $\|A\|_2\approx 1$ (this can be enforced using an inexpensive norm estimator $\alpha\approx \|A\|_2$, followed by a scaling $A\leftarrow A/\alpha$), the idea is simply to compute 
$r(A):=Ur(\Sigma)V^*$,
where $r$ is a rational function that maps the interval $[s,1]$ to $1$, to machine precision. 
The upshot is that if $s\gg \sigma_{\min}(A)$, 
then $r$ is allowed to be of much lower degree than would be needed for 
computing $U_p$. See Fig.~\ref{fig:rxsvd} for an illustration. 

It is worth noting that the degree of $r$ and the number of QDWH iterations depend on the user-defined value $s$, unlike \texttt{QDWHpartial-EIG} (for which $s=0.2$ is a ``fixed" choice). Accordingly, if $s\ll 1$ (say $s<10^{-3}$), for the first QDWH iteration it is advisable to use the QR-based implementation in Eq.~\eqref{practical-polar-qr} rather than Eq.~\eqref{practical-polar-chol} to avoid instabilities. 

We also note that we use $|\ell_k-1|<O(u)$ as the stopping criterion for QDWH. This is because the standard condition, which requires convergence of $X_k$, is not necessarily satisfied when $[s,1]$ has been mapped to $1\pm O(u)$, because the singular values $\sigma_i(\tilde A)<s$ have not converged and lie somewhere in $[0,1]$. (The situation was the same in \texttt{QDWHpartial-EIG}, but there it was simpler as the iteration number is always three.)

Once $r(A)$ is computed, one can find the desired column space $U_1$ 
(leading columns of $U$ corresponding to singular values of $A$ larger than $s$) by finding the null 
space of $I-r(A)^*r(A)$, which we do as before using the QR factorization (with randomization if needed) and looking for the diagonal entries of $R$. 
It is worth noting that computing the matrix $r(A)^*r(A)$ may seem ill-advised, as it squares the condition number. This is not an issue here, as the quantity of interest is the singular subspace corresponding to the largest singular values of $r(A)$ (which are $\approx 1$).

Algorithm~\ref{algo:partial-svd} gives a pseudocode of the overall algorithm. 
\begin{algorithm}
	\footnotesize{
	\begin{algorithmic}[1] 
\State Estimate $\alpha \approx \|A\|_2$ with a norm estimator.
\State Apply QDWH to 
$\tilde A: = A/\alpha$ with $\ell_0=s$ until $|\ell_k-1|<O(u)$ to obtain $r(\tilde A)$.
	\State Calculate $[Q\; R] = QR(I-r(\tilde A)^*r(\tilde A))$.
	\State Find the index $ind = \min(find(abs(diag(R))<tol=0.01))$.
	\State Extract $Q_2 = Q(:,ind:end)$.
 	\State Projection: Compute $AQ_2$ and its SVD $AQ_2=\tilde U\tilde\Sigma \tilde V^*$. 
  \State Extract the singular triplets with singular values larger than $s\|A\|_2$, call them $\tilde U_1,\tilde\Sigma_1, \tilde V_1$. 
 	\State Let $U_1:=\tilde U_1,\Sigma_1:=\tilde\Sigma_1$ and $V_1:=Q\tilde V_1$. 
	\end{algorithmic}
	}
\caption{\texttt{QDWHpartial-SVD}. Given $A\in\mathbb{C}^{m\times n}$ ($m\geq n$) and threshold $s>0$, compute the singular values larger than $s\|A\|_2$ and the corresponding singular vectors. 
} 
\label{algo:partial-svd}
\end{algorithm}

The outputs of \texttt{QDWHpartial-SVD} are $U_1,\Sigma_1,V_1^*$ such that $U_1\Sigma_1V_1^*$ is the truncated SVD of $A$, truncated at the first singular value smaller than $s$. 

The computational savings comes from the fact that $AQ_2$ is much thinner than $A$; the number of columns of $AQ_2$ is slightly more than the number of singular values of $A$ larger than $s$. 

Of course, if $m<n$ one can simply apply the algorithm to $A^*$. 


\section{Implementation details}
\label{sec:impl}
The main basic blocks of \texttt{QDWHPartial-EIG} in Algorithm~\ref{algo:partial-eig}
correspond to a subset of \texttt{QDWHPartial-SVD} operations, as
shown in in Algorithm~\ref{algo:partial-svd}.
Therefore, we only provide implementation details of \texttt{QDWHpartial-SVD}.
Algorithm~\ref{algo:qdwh-partial-svd} describes the pseudo-code of the distributed-memory implementation
of \texttt{QDWHpartial-SVD} based on \texttt{ScaLAPACK}~\cite{ScaLAPACK}.

Following the 2D Block-Cyclic Data Distribution (2D-BCDD) used in \texttt{ScaLAPACK}, we define 
the MPI process grid configuration as $P \times Q$. Each data structure owns a handle or a descriptor
that expresses how the data structure is distributed following the 2D-BCCD.
\texttt{ScaLAPACK} relies on the Basic Linear Algebra Communication Subprograms (BLACS) library, which is in charge
of performing data movements during the matrix computations through the traditional MPI.
\texttt{ScaLAPACK} relies on block algorithms, which can be expressed by two successive computational stages: 
the panel factorization and the update of the trailing submatrix. While the former is memory-bound, 
and typically sequential and may not benefit from having many processors participating,
the latter is rich in compute-bound operations. This is where most of \texttt{ScaLAPACK} 
dense linear algebra operations extract parallel performance
by means of calls to Level-3 BLAS, as implemented in the Parallel BLAS (PBLAS) layer.
The blocking size referred as \emph{nb} is an internal tuning parameter that 
trades-off the degree of parallelism and the performance of the computational kernels.
Moreover, the number of processors should be properly calibrated into a rectangular shape
with $Q > 1.5 P$ to carry on, in parallel, the update of the trailing submatrix.

As shown in Algorithm~\ref{algo:qdwh-partial-svd}, the \texttt{QDWHpartial-SVD} code
is mostly composed of conventional dense linear algebra
matrix kernels rich in compute-intensive Level-3 BLAS operations that are capable of achieving a 
decent percentage of the system's theoretical peak performance. Since these matrix kernels are 
widely available in vendor optimized numerical libraries, porting to various hardware architectures 
should not be cumbersome.
The code is written in double precision arithmetics, and can be extended to other 
precisions for a broader application coverage. 

\begin{algorithm}
\scriptsize{
\begin{algorithmic}[1] 
\Statex \textbf{/* Set block size and initiate CBLACS context */}
\State $Cblacs\_get(0, 0, ictxt)$
\State $Cblacs\_gridmap(ictxt, imap, P, P, Q)$
\State $Cblacs\_gridinfo(ictxt, P, Q, myrow, mycol)$
\State \textbf{/* Initialize data structures using the 2D-BCDD descinit() */}
\State $descinit(nb, nb, A, descA); Fill\_in(A, descA)$
\Statex \textbf{/* Estimate the two-norm of the matrix */}
\State $\alpha = pdgenm2(A)$
\State $pdlascl(\alpha, 1., A)$
\Statex \textbf{/* Computing the polar factor $U_{p}$ of the matrix $A$ using QDWH */}
\State $k = 1,\; Li = threshold (s), \; conv = 100$

\While {$(\vert Li - 1 \vert \ge 5eps)$}
\State $L2= Li^{2}, \; dd = \sqrt[3]{(4(1-L2)/L2^{2})}$   
\State $sqd = \sqrt{1+dd}$
\State $a1 = sqd + \sqrt{8 - 4 \times dd + 8(2-L2)/(L2 \times sqd))}/2$
\State $a = real(a1); \; b = (a-1)^{2}/4; \; c = a + b -1$
\State $Li = Li (a + b \times L2)/(1 + cL2)$
\State pdlacpy(U, U1) 
\State \textbf{/* Compute $U_{k}$ from $U_{k-1}$ */}

\State $pdlaset(Z, 0., 1.)$ 
\State $pdgemm(U^{\top}, U, Z)$ 
\State $pdgeadd(U, B)$ 
\State $pdposv(Z, B)$ 
\State $pdgeadd(B, U)$ 

\State $pdgeadd(U, U1)$ 
\State $pdlange(U1, conv)$ 
\State $k = k + 1$
\EndWhile
\Statex \textbf{/* $U_{p}$ contains the isolated subspectrum of interests */}

\State $pdlaset(0.0, 1.0, B)$
\State $pdgemm(U_{p}^{\top}, U_{p}, B)$
\State $pdgeqrf(B)$
\State $ind = min(find(abs(diag(B)) < tol = 0.01))$
\State $pdorgqr(B, Q)$
\Statex \textbf{/* size$(\tilde{A}) = N - ind$ */}
\State $Q_{2} = Q({:},ind{:}end)$
\State $pdgemm(A, Q_{2}, \tilde{A})$
\Statex \textbf{/* Calculate the SVD on the reduced problem */}
\State $pdgesvd(\tilde{A}, \tilde{U}, \tilde{\Sigma}, \tilde{V})$
\Statex \textbf{/* $\tilde{U_{1}}, \tilde{\Sigma_{1}}, \tilde{V_{1}}$ are the singular triplets with singular values larger than the $s \times \alpha$ */}
\State $U = \tilde{U_1}$
\State $\Sigma = \tilde{\Sigma_1}$
\State $pdgemm(\tilde{V_1}, Q_{2}^T, V)$

\end{algorithmic}
}
\caption{Pseudo-code of the \texttt{QDWHpartial-SVD} using \texttt{ScaLAPACK}.}
\label{algo:qdwh-partial-svd}
\end{algorithm}

\section{Operation counts}
\label{sec:complexity}
Table~\ref{table:complexity} reports the operation counts of various symmetric EIG and SVD
solvers  on square matrices of size $n$: the \texttt{PDSYEVD} / \texttt{PDGESVD} and
\texttt{QDWH-EIG} / \texttt{QDWH-SVD} routines for computing the full spectrum and the
\texttt{QDWHpartial-EIG} / \texttt{QDWHpartial-SVD} routines for computing a subset of the spectrum. 
We refer the reader to~\cite{Yuji2013} for further details on the costs of the standard and QDWH-based EIG / SVD solvers.

The operation counts of \texttt{QDWHpartial-EIG} and \texttt{QDWHpartial-SVD} 
depends on the number of Cholesky-based QDWH iterations (typically two or three) and, 
the \texttt{QR}, \texttt{GEMM} and \texttt{SYRK} to form the reduced problem matrix of size $N_{s} \ge s$, with 
$s$ the size of the partial spectrum of interest. The actual full EIG and SVD 
occurs now only on the reduced problem matrix of size $N_{s}$.
Assuming $N_{s}\ll N$ and three Cholesky-based iterations to get the polar factor from QDWH, 
the total number of operations is up to $14N^3$ and $24N^3$ for \texttt{QDWHpartial-EIG} ($it_{Chol}=3$)
and \texttt{QDWHpartial-SVD} ($it_{QR}=1$ and $it_{Chol}=3$),
respectively.
\begin{table}[htbp]
    \begin{center}
	\scriptsize{
    \begin{tabular}{ c | c }
 EIG and SVD variants&  Cost  \\
 \hline\hline
  & \\
 Standard full EIG & $9N^{3}$ \\
  & \\
 Full QDWH-EIG &   $(17+\frac{4}{9})N^{3} \le$ $\cdots$  $\le (52+\frac{1}{9})N^{3}$  \\
  & \\
 \multirow{3}{*}{QDWHpartial-EIG} & QDWH: (4+1/3)$N^{3}$ x $\#it_{Chol}$  \\
                                  & QR + SYRK + GEMM: $4/3N^{3}$ + $N_{s}N^2$ + $2N_{s}^2N$\\
                                  & EIG: $9N_{s}^{3}$\\
  & \\
 \hline\hline
  & \\
 Standard full SVD &  $17N^{3}$   \\
  & \\
 Full QDWH-SVD &   $20N^{3} \le$ $\cdots$  $\le (50  + \frac{1}{3})N^{3}$  \\
  & \\
  \multirow{3}{*}{QDWHpartial-SVD} & QDWH: (8+2/3)$N^{3}$ x $\#it_{QR}$ + (4+1/3)$N^{3}$ x $\#it_{Chol}$ \\
                                   & QR + SYRK + $2\times$GEMM: $4/3N^{3}$ + $N^3$ + $4NN_{s}^2$\\
                                   & SVD: $17N_{s}^{3}$\\
  & \\
 \hline\hline
    \end{tabular}
	}
    \end{center}
	\caption{Operation counts for various symmetric EIG and SVD algorithms.}
	\label{table:complexity}
\end{table}


\section{Numerical accuracy}
\label{sec:acc}
The numerical accuracy of the QDWH-based algorithms to compute the polar decomposition,
the eigenvalue decomposition (\texttt{QDWH-EIG}) and singular value decomposition (\texttt{QDWH-SVD}) 
have been verified in~\cite{Yuji2013}. The robustness of their high performance implementations
has been studied on shared-memory systems~\cite{Sukkari2016}
and on distributed-memory systems~\cite{Sukkari2016bis,Sukkari2017bis}.
In this Section, we present the numerical robustness of the 
\texttt{QDWHpartial-EIG} and \texttt{QDWHpartial-SVD} implementations on distributed-memory
system. 
\subsection{Environment Settings}
\label{sec:env_set}
We run our experiments on a Cray XC40 system codenamed \emph{Shaheen-2}
installed at the KAUST Supercomputing Laboratory (KSL), with the 
Cray Aries network interconnect, which implements a Dragonfly 
network topology.
The system has $6174$ compute nodes, each with two-socket 16-core Intel Haswell 
running at $2.3$GHz and $128$GB of DDR3 main memory. 
The Haswell nodes on \emph{Shaheen-2} have a theoretical peak performance of
approximately 1.18 TFlops/s.
Furthermore, ``hugepages" are employed to improve memory accesses.
The work load managers on \emph{Shaheen-2} is native
SLRUM. 
We use the Intel compiler v15.0.2.164.
We rely on the \texttt{ScaLAPACK} implementation from the high performance Cray LibSci
numerical library, which depends on the MPI programming model for inter-node communications.\\
All runs for a given process configuration have been submitted in the same 
job submission script to reduce the impact from the system jitter.
All the experiments are performed using IEEE double-precision arithmetic.
\subsection{Synthetic Matrices}
\label{subsec:synth_matrix}
The dense synthetic matrices $A \in \mathbb{R}^{N \times N}$ are generated using the \texttt{ScaLAPACK} 
routine \texttt{PDLATMS} $A = Q_{1}D\mathnormal{Q_{2}}^\top$ with setting mode = 0. 
For the symmetric EIG solvers testing, the matrices are generated with an equispaced eigenvalues as follows:
\[
  \begin{tabular}{cccccc}
  \multirow{2}{*}{D =} \rdelim\{{2}{10pt}[] {}& $D[i]$ =& $-k \times randn[i]$, & $i \le k $ & \rdelim\}{2}{20pt}[,]\\
  & $D[i]$ =& $N-k \times randn[i]$, & $i \ge k $ & &\\
  \end{tabular}
\]
where $k$ is the number of the negative eigenvalues.
For the SVD solvers testing, the distribution of the singular values of the generated matrices follows 
a geometrical series: $D[i] = (0.5) ^ {\frac{i}{N}*100}$ 
We compute then orthogonal matrices $Q_{1}$ and $Q_{2}$ generated by 
calculating the $QR$ factorization of arbitrary matrices to form the SVD,
while $Q_{1} = Q_{2}$ for the symmetric EIG solvers. The performance of the 
matrix generation step may be expensive and can be improved but this is 
beyond the scope of this paper.
\subsection{Norm Definitions}
\label{subsec:norm}
For a given general matrix $A \in \mathbb{R}^{N \times N}$, let $\Sigma = diag(\sigma_{1}, \sigma_{2}, ..., \sigma_{k})$
be the $k$ computed singular values, and $U$ and $V$
be the corresponding $k$ computed left and right singular vectors.
The norm $\Vert\; .\;\Vert_{F}$ denotes the Frobenius norm. The accuracy assessment of the 
partial computation of the SVD are based on the 
following metrics: 
\begin{equation}
\frac{\Vert I-UU^{\top} \Vert_{F}}{n} \; \hbox{\rm and} \; \frac{\Vert I-VV^{\top} \Vert_{F}}{n},
\label{eq:orth}
\end{equation}
for the orthogonality of the left and right $k$ computed singular vectors $U$ and $V$,
respectively, and 
\begin{equation}
\frac{\Vert \Sigma-\Delta \Vert_{F}}{\Vert \Delta \Vert_{F}}, 
\label{eq:sv}
\end{equation}
for the accuracy of the $k$ computed singular values $\Sigma$, where $\Delta$ 
is the $k$ exact singular values (analytically known), 
and
\begin{equation}
\max_{i}(   \Vert AU(:,i)  -  \sigma_{i}V(:,i)    \Vert_{F}    ) 
\; \text{and} \;
\max_{i}(   \Vert AV(:,i)  -  \sigma_{i}U(:,i)    \Vert_{F}    )
\label{eq:allsv}
\end{equation}
for the accuracy of the left and right singular value decomposition, respectively.
Similarly, for a symmetric matrix $A \in \mathbb{R}^{n \times n}$, 
the accuracy of the computed $k$ negative eigenvalues, 
the orthogonality of their corresponding eigenvectors and the overall residual
can be accordingly measured using~\ref{eq:orth},~\ref{eq:sv} and~\ref{eq:allsv}, respectively.
\subsection{Accuracy Assessments of EIG/SVD Solvers}
\label{subsec:acc-svd}
\begin{figure}[htpb]
      \centering
      \begin{minipage}[htpb]{.46\textwidth}
          \subfigure[Accuracy of singular values.]{
			 \includegraphics[width=.92\linewidth]{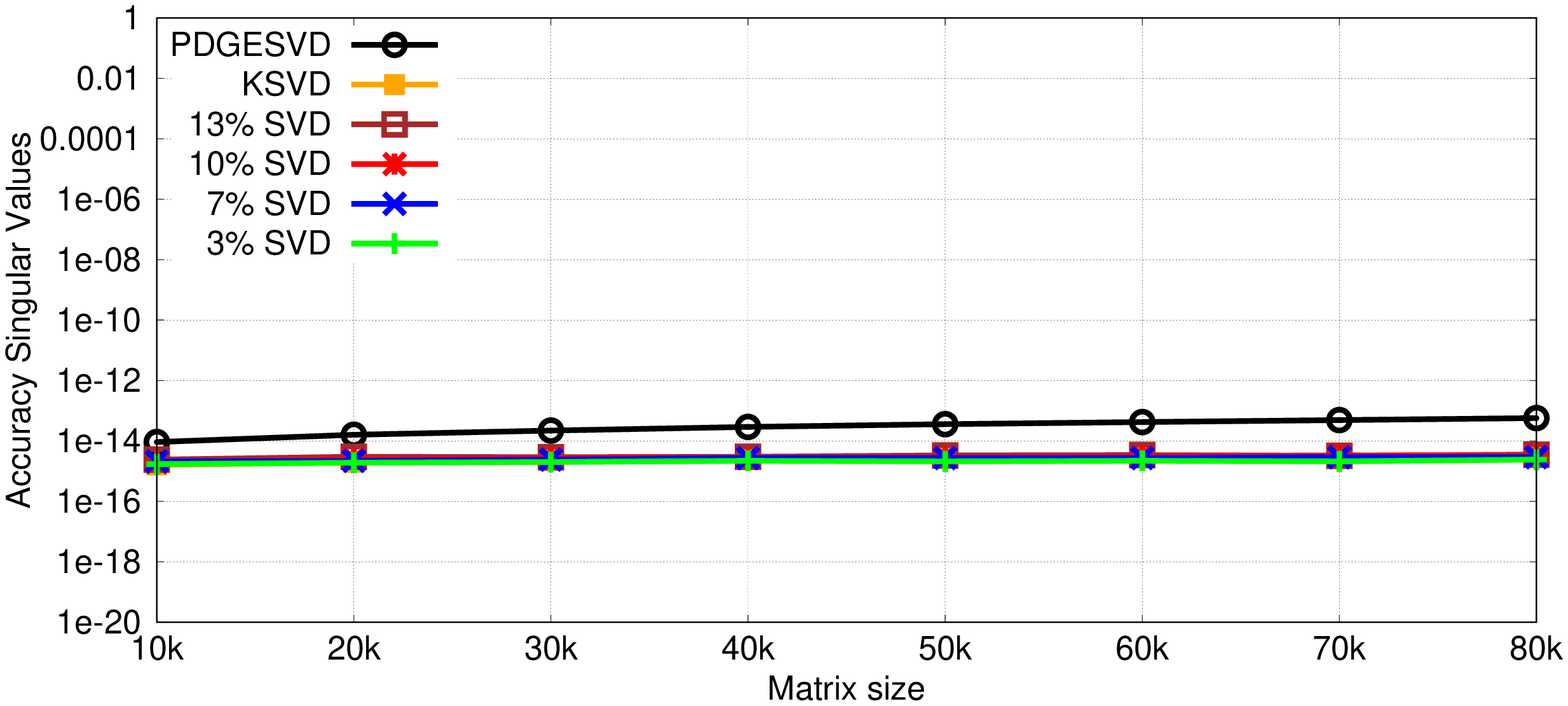}}
          \subfigure[Orthogonality of left/right singular vectors.]{
			 \includegraphics[width=.92\linewidth]{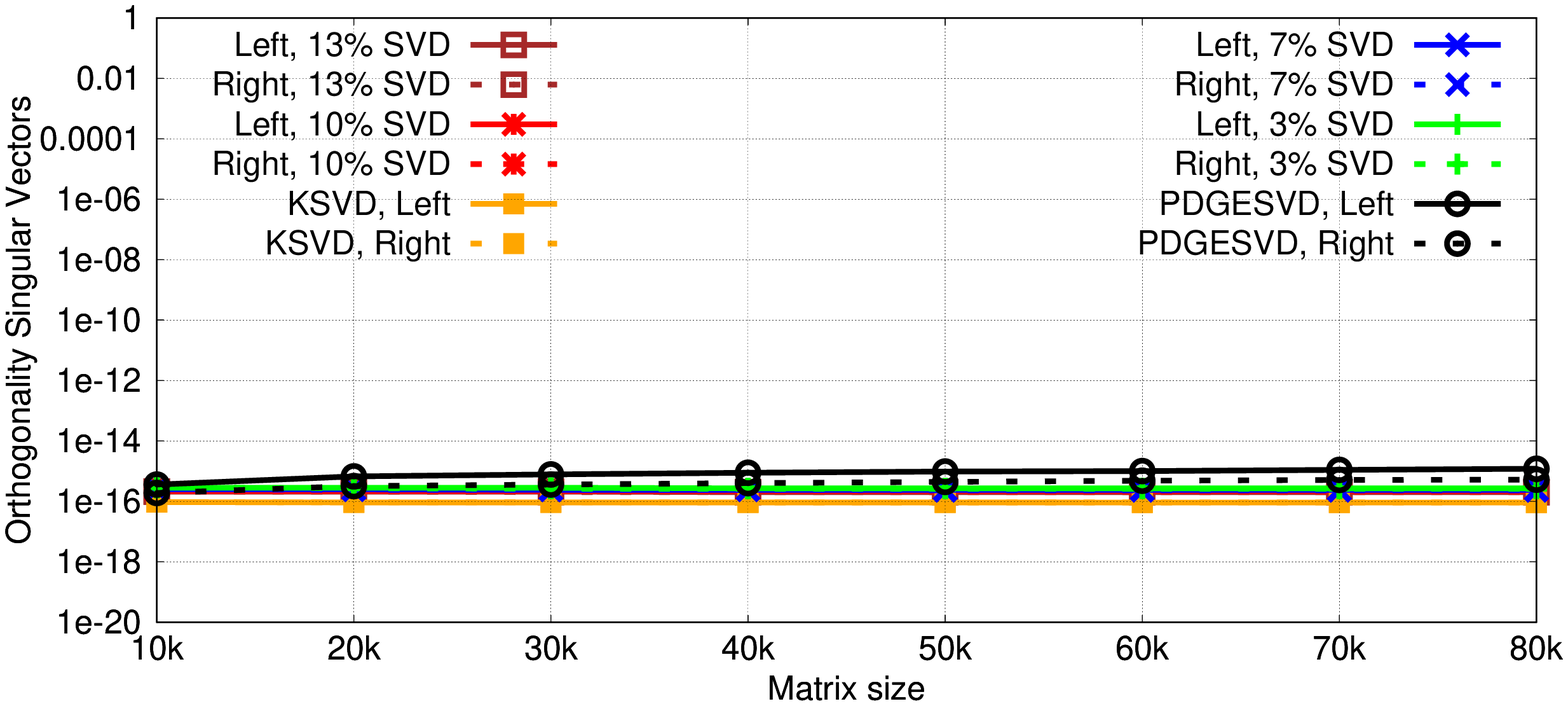}}
          \subfigure[Backward stability.]{
			 \includegraphics[width=.92\linewidth]{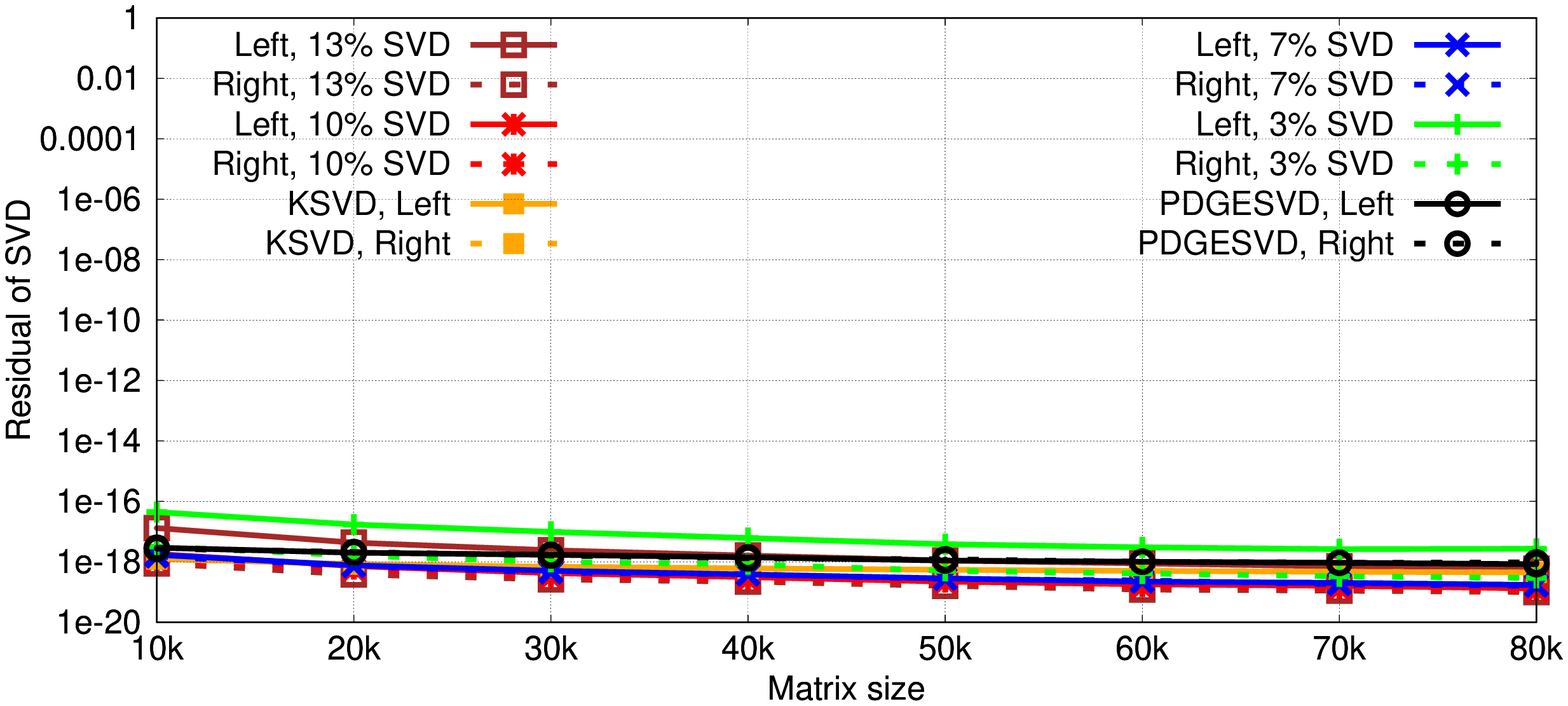}}
      \end{minipage}
      \begin{minipage}[htpb]{.46\textwidth}
          \subfigure[Accuracy of eigen values.]{
			 \includegraphics[width=.92\linewidth]{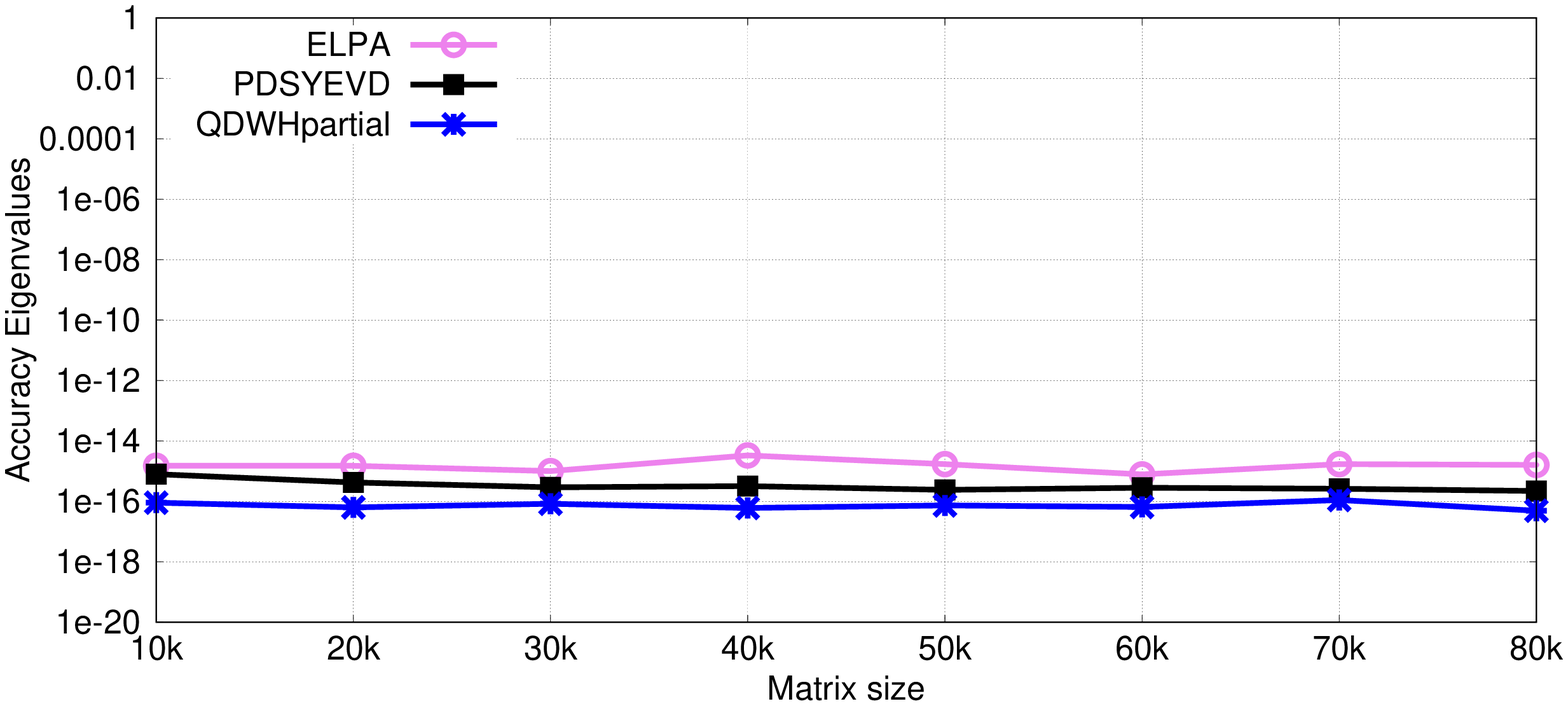}}
          \subfigure[Orthogonality of eigen vectors.]{
			 \includegraphics[width=.92\linewidth]{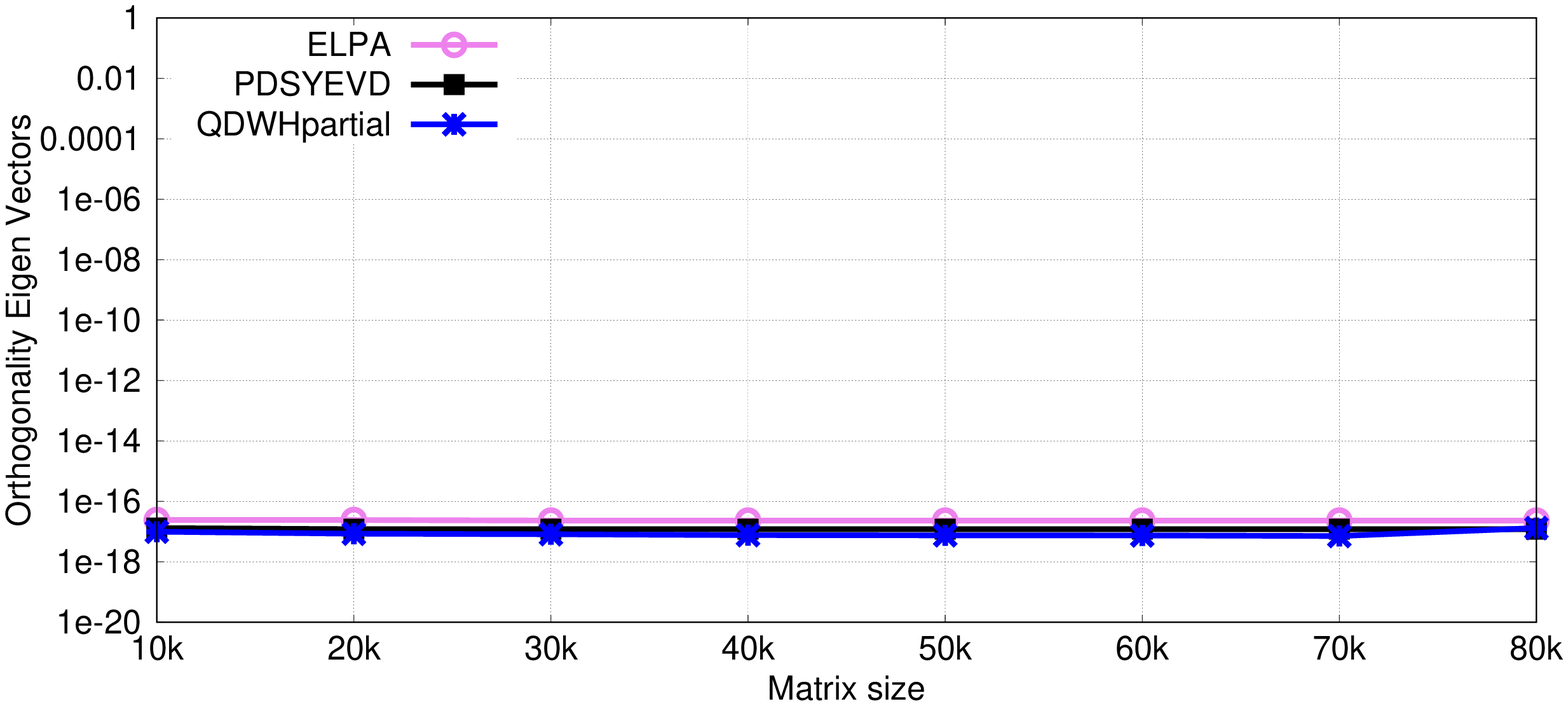}}
          \subfigure[Backward stability.]{
			 \includegraphics[width=.92\linewidth]{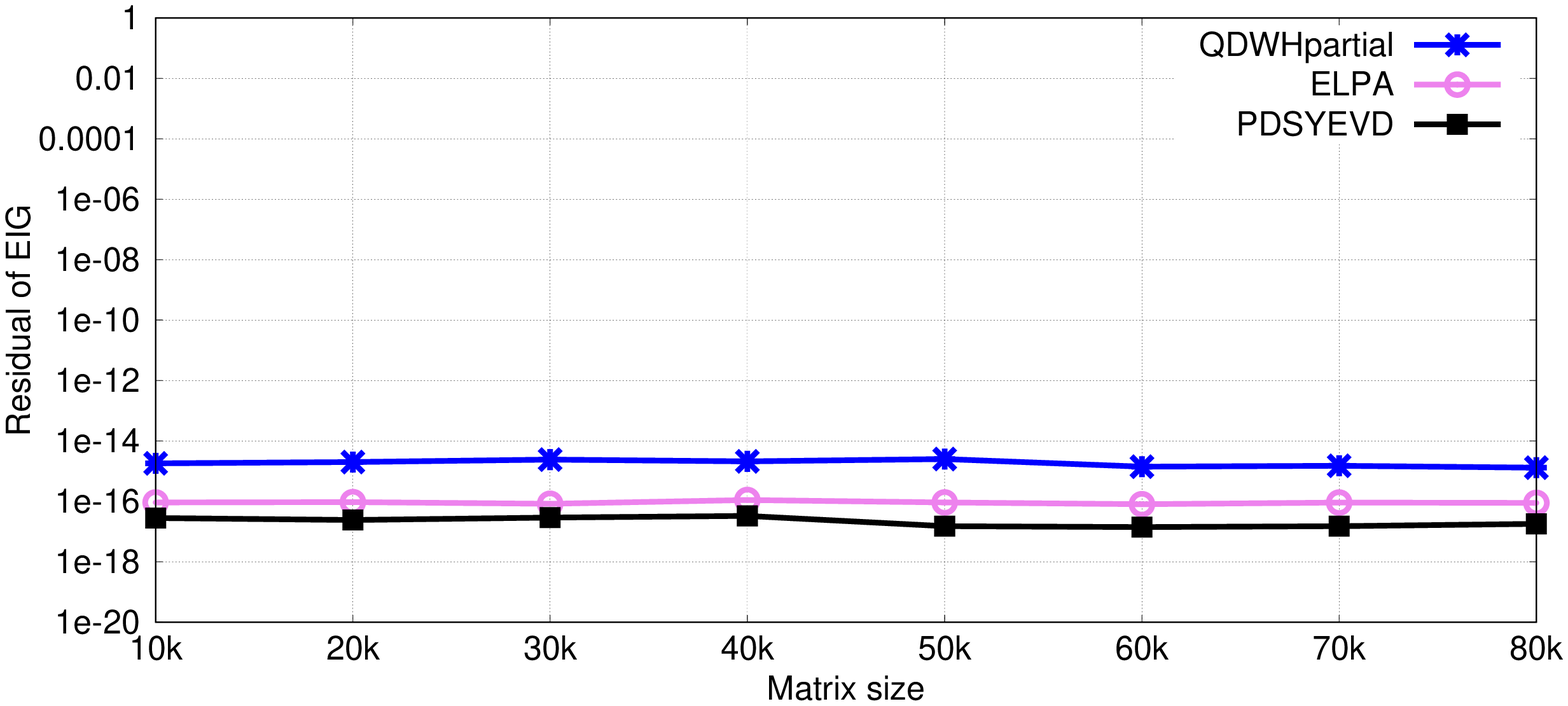}}
      \end{minipage}
      \caption{Assessing the numerical accuracy/robustness using $16 \times 36$ grid topology: (a-b-c) for SVD solvers and (d-e-f) for EIG solvers}
      \label{fig:num-18}
\end{figure}

This section highlights the numerical robustness of
\texttt{QDWHpartial-SVD} and \texttt{QDWHpartial-EIG} implementations
Fig.~\ref{fig:num-18} (a, b, c) shows the numerical accuracy
of the computed singular values (Equation~\ref{eq:sv}), 
the orthogonality of their corresponding singular vectors (Equation~\ref{eq:orth}) 
and the right/left residual of the computed 
SVD (Equation~\ref{eq:allsv}) 
on a $16 \times 36$ grid configuration (similar accuracy results for larger grid sizes
 $32 \times 72$, $64 \times 144$ and $128 \times 288$)
using synthetic ill-conditioned matrices. 
Herein, we compare the accuracy of three implementations of the SVD solvers: 
\texttt{QDWHpartial-SVD} (setting different threshold $s$) against
\texttt{PDGESVD} from \texttt{ScaLAPACK} and from \texttt{KSVD}\footnote{Available at \url{https://github.com/ecrc/ksvd}}~\cite{Sukkari2017bis}.
The \texttt{QDWHpartial-SVD} is capable to extract only the singular 
values/vectors of interest within the user-defined threshold $s$.
This threshold $s$ can be tuned with \emph{a priori} knowledge on the singular value distribution (e.g., 
globally low-rank structure).
This tunable parameter can directly influence the number of the computed singular values/vectors.
For instance, in Fig.~\ref{fig:num-18} (a, b, c), we study the accuracy 
for $s = 0.1, 0.01, 0.001, 0.0001$ that translates into the percentages
3\%, 7\%, 10\%, 13\% of the computed singular values/vectors, 
respectively, and as a result affect the performance of \texttt{QDWHpartial-SVD}. 
It is noteworthy that the \texttt{ScaLAPACK} \emph{PDGESVD} 
computes first the whole SVD, then the requested singular values/vectors are 
filtered out using the threshold parameter $s$.

Figures~\ref{fig:num-18} (d, e, f)
shows the numerical accuracy of three different EIG solvers to compute 10\% of the negative eigenvalues: 
\texttt{ELPA} divide-and-conquer routine, \texttt{ScaLAPACK} \emph{PDSYEVD} and \texttt{QDWHpartial-EIG}. 
\texttt{ELPA} and QDWHpartial-EIG are capable of computing a fraction of 
the negative eigenspectrum. The \texttt{ScaLAPACK} \emph{PDSYEVD} calculates first 
the entire eigenspectrum and then only the 10\% of the negative eigenspectrum are selected.

These extensive numerical tests in Fig.~\ref{fig:num-18} demonstrate the numerical 
robustness of \texttt{QDWHpartial-SVD} and \texttt{QDWHpartial-EIG} to provide satisfactory accuracy up to 
the machine precision for double precision computations across all studied matrix sizes.

\section{Performance results}
\label{sec:results}
\begin{figure}[htpb]
      \centering
      \begin{minipage}[htpb]{.46\textwidth}
          \subfigure[P=16 and Q=36.]{
                         \includegraphics[width=.92\linewidth]{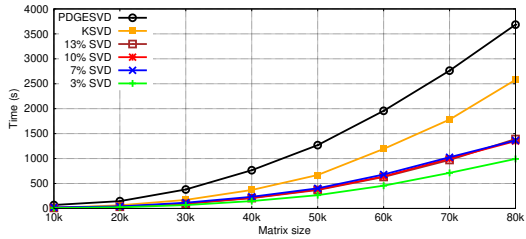}}
          \subfigure[P=32 and Q=72.]{
                         \includegraphics[width=.92\linewidth]{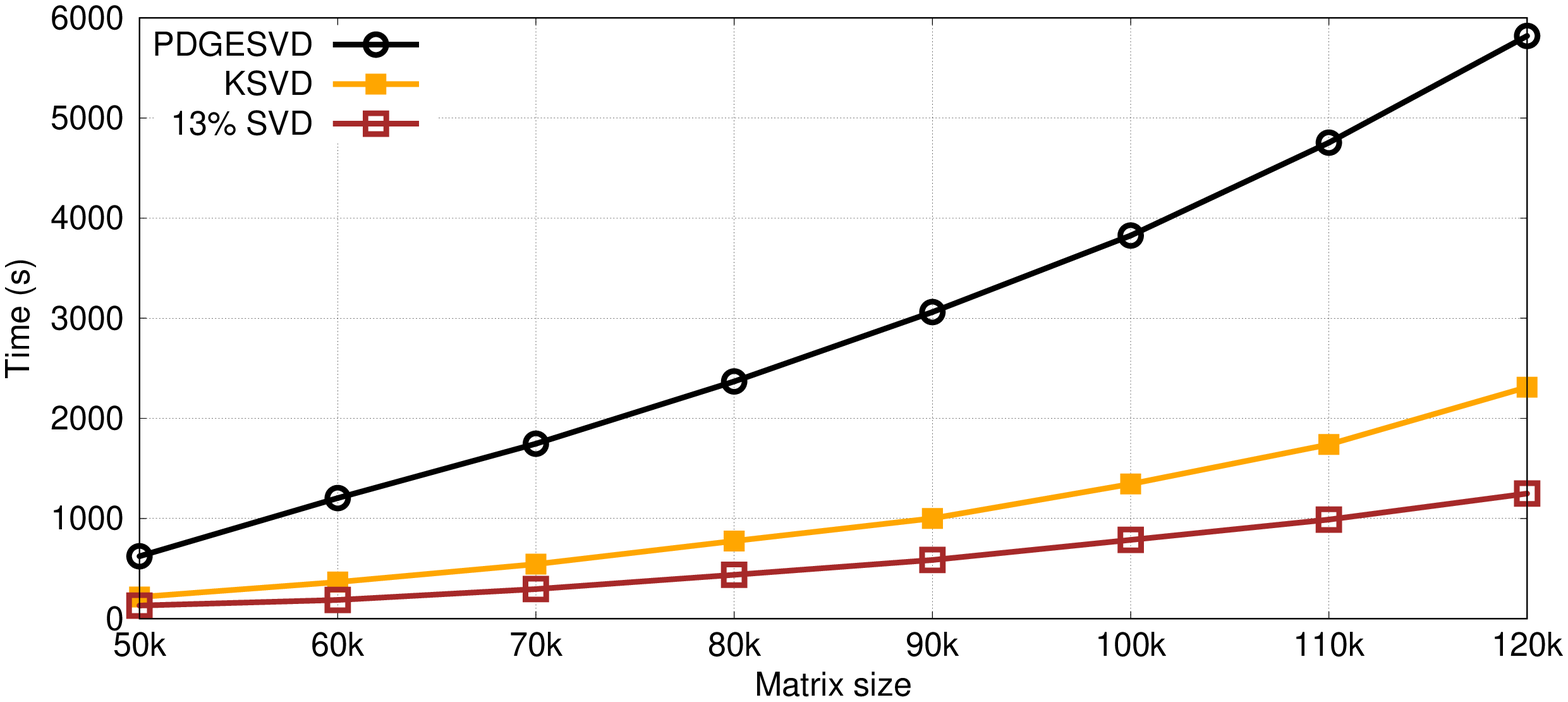}}
      \end{minipage}
      \begin{minipage}[htpb]{.46\textwidth}
          \subfigure[P=64 and Q=144.]{
                         \includegraphics[width=.92\linewidth]{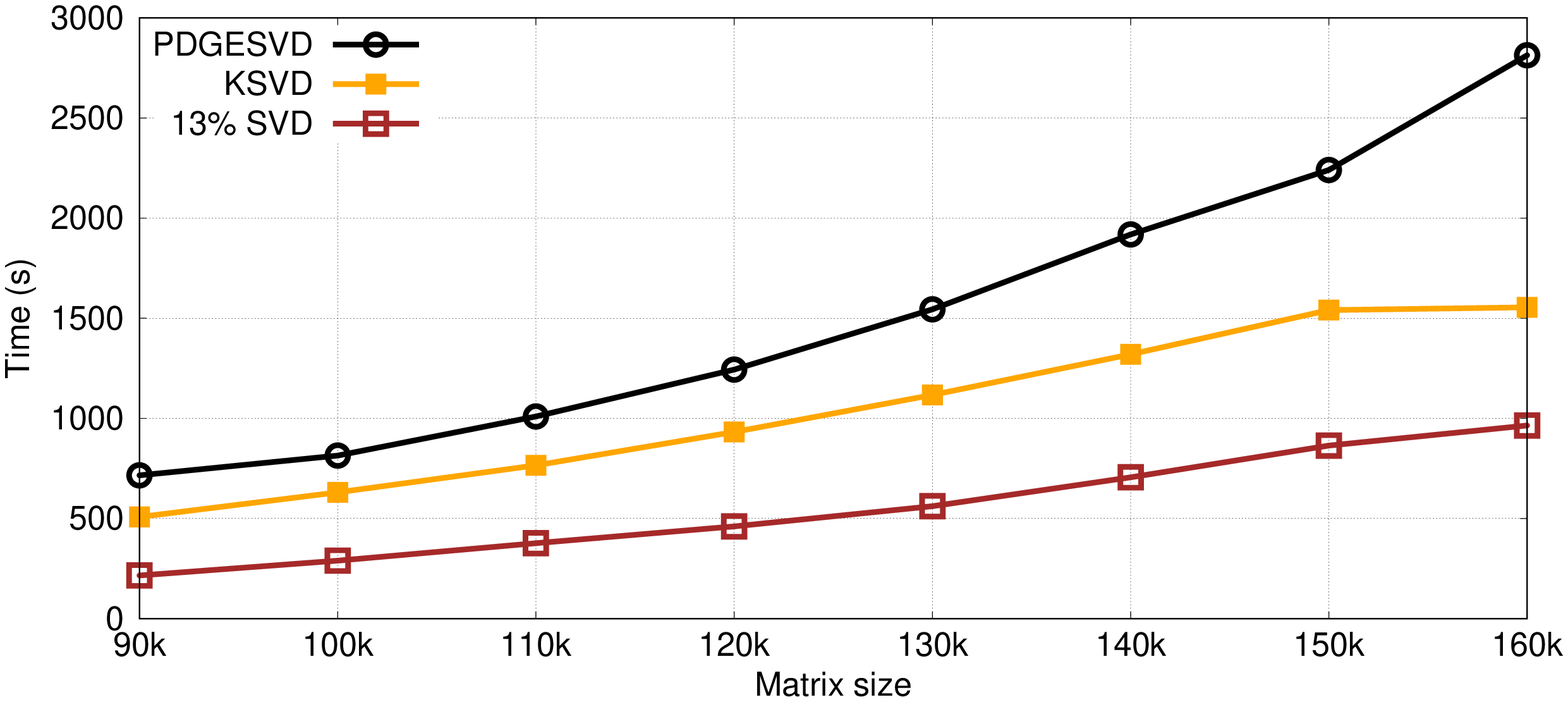}}
          \subfigure[P=128 and Q=288.]{
                         \includegraphics[width=.92\linewidth]{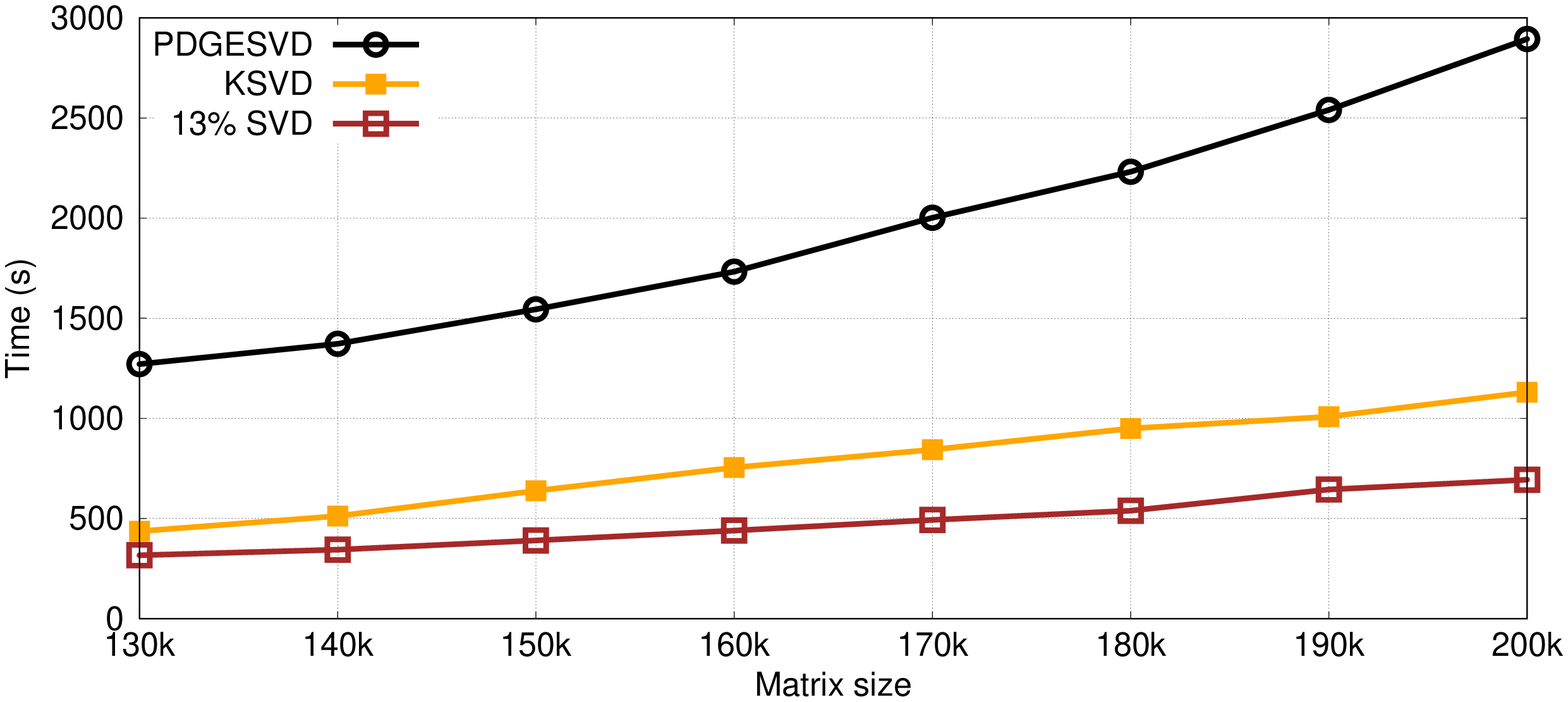}}
      \end{minipage}
      \caption{Performance comparisons in seconds of SVD solvers for various grid topologies.}
      \label{fig:svd-comparison-time}
\end{figure}

Figure~\ref{fig:svd-comparison-time} highlights the performance comparisons of
\texttt{QDWHpartial-SVD} against the two other SVD solvers,
\texttt{ScaLAPACK} \emph{PDGESVD} and \texttt{KSVD}~\cite{Sukkari2017bis} across various matrix 
sizes and process grid configurations. 
When only 13\% of the spectrum is needed, 
\texttt{QDWHpartial-SVD} achieves performance superiority as we increase the matrix sizes
up to [$6.3X$, $2X$] on $16 \times 36$, [$6X$, $2X$] on $32 \times 72$,
[$3.3X$, $2.3X$] on $64 \times 144$ and [$4X$, $1.8X$] on $128 \times 288$ grid topologies against
[\texttt{ScaLAPACK} \emph{PDGESVD}, \texttt{KSVD}], respectively.
Moreover, on $16 \times 36$ grid configuration, \texttt{QDWHpartial-SVD} achieves similar performance
when only 13\%-10\%-7\% of the spectrum is needed, since the QDWH-based polar decomposition
is the most time consuming step. We observe a slightly faster time to solution though, when only 3\%
of the spectrum is calculated, since the size of the reduced problem maybe relatively smaller than
the aforementioned partial spectrums.

\begin{figure}[htpb]
      \centering
      \begin{minipage}[htpb]{.46\textwidth}
          \subfigure[P=16 and Q=36.]{
                         \includegraphics[width=.92\linewidth]{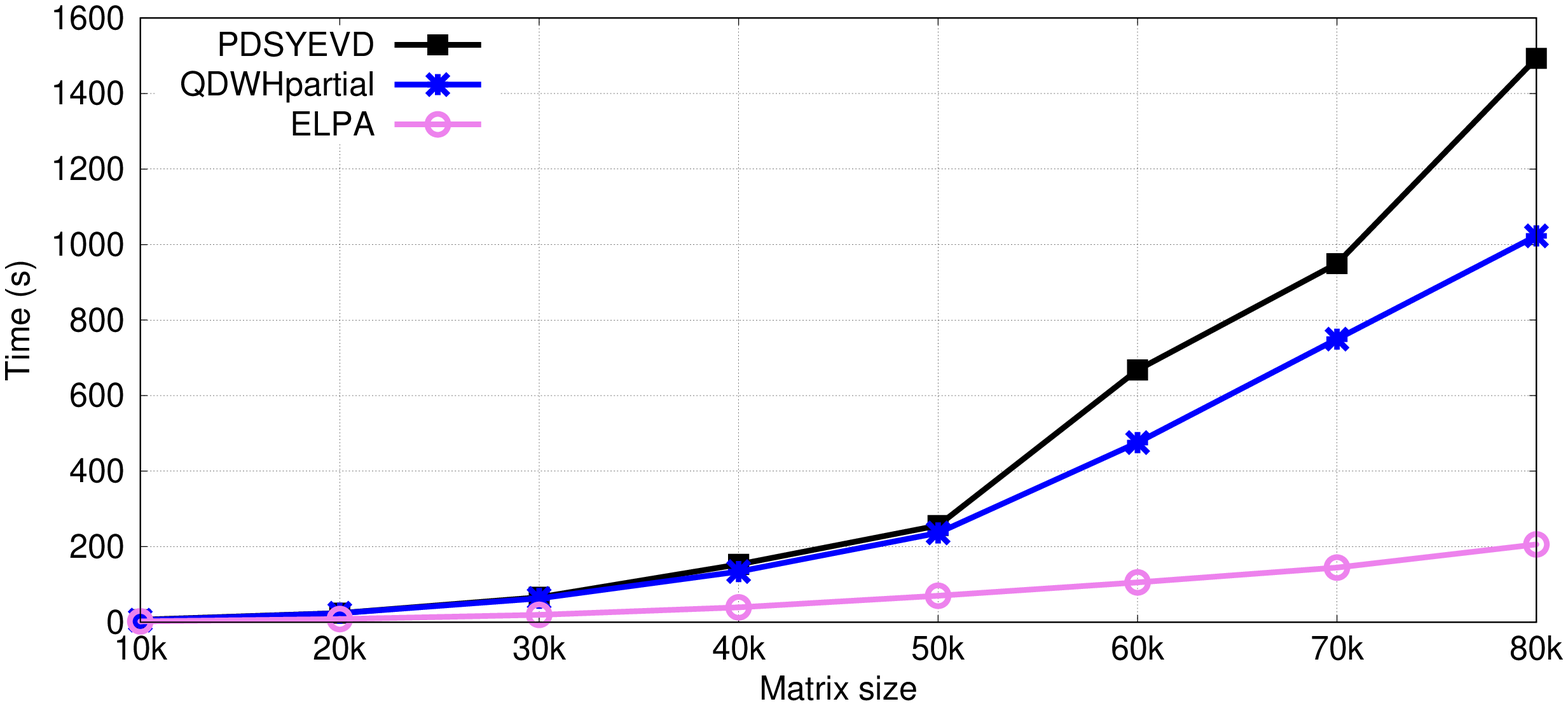}}
          \subfigure[P=32 and Q=72.]{
                         \includegraphics[width=.92\linewidth]{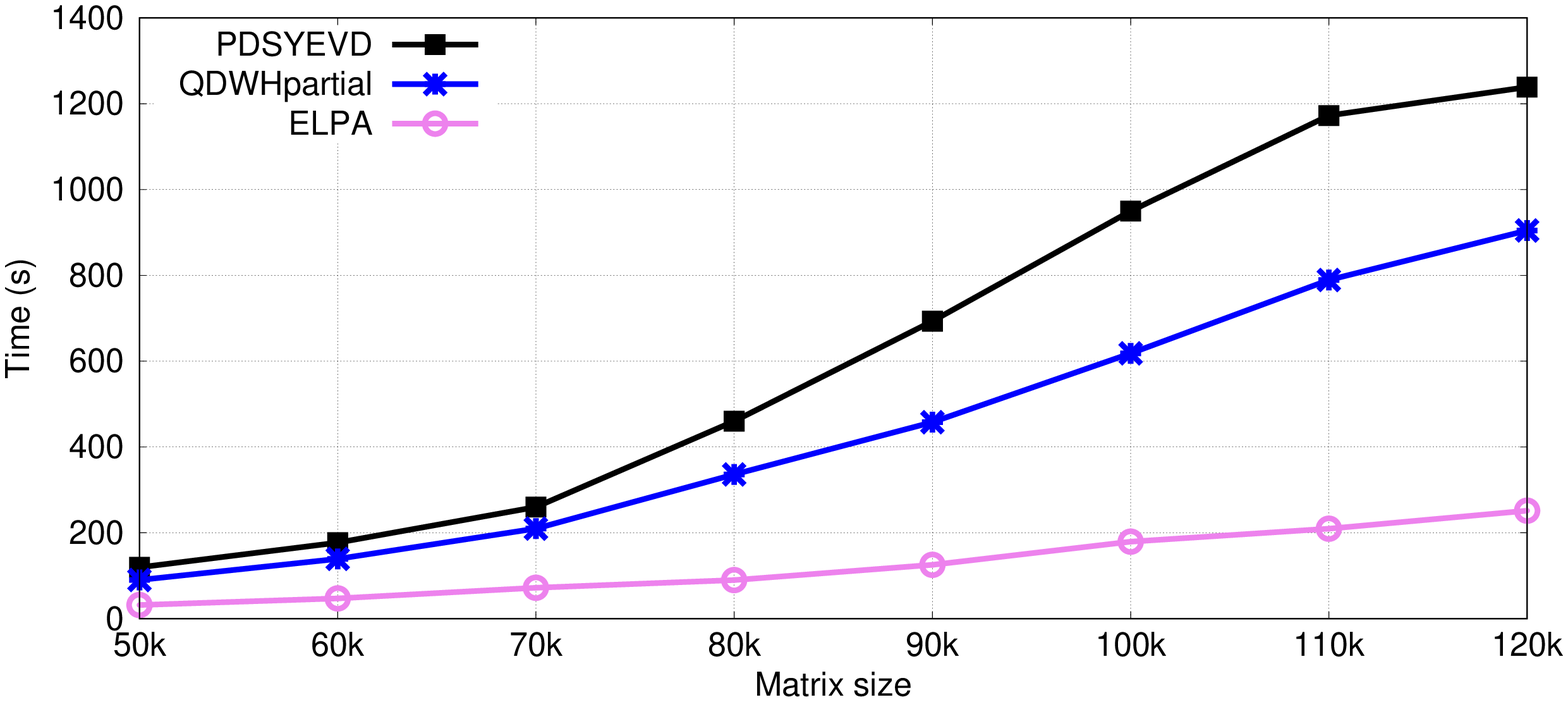}}
      \end{minipage}
      \begin{minipage}[htpb]{.46\textwidth}
          \subfigure[P=64 and Q=144.]{
                         \includegraphics[width=.92\linewidth]{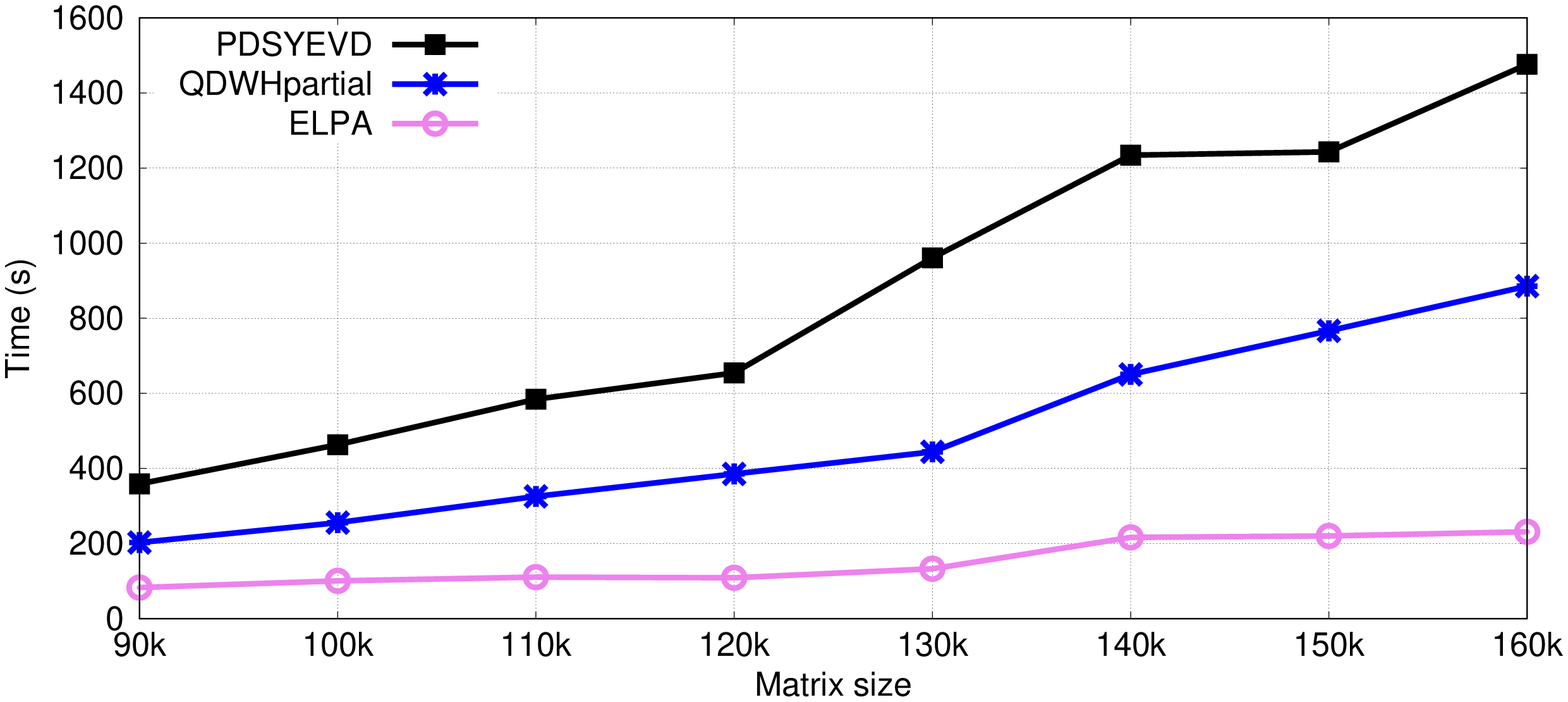}}
          \subfigure[P=128 and Q=288.]{
                         \includegraphics[width=.92\linewidth]{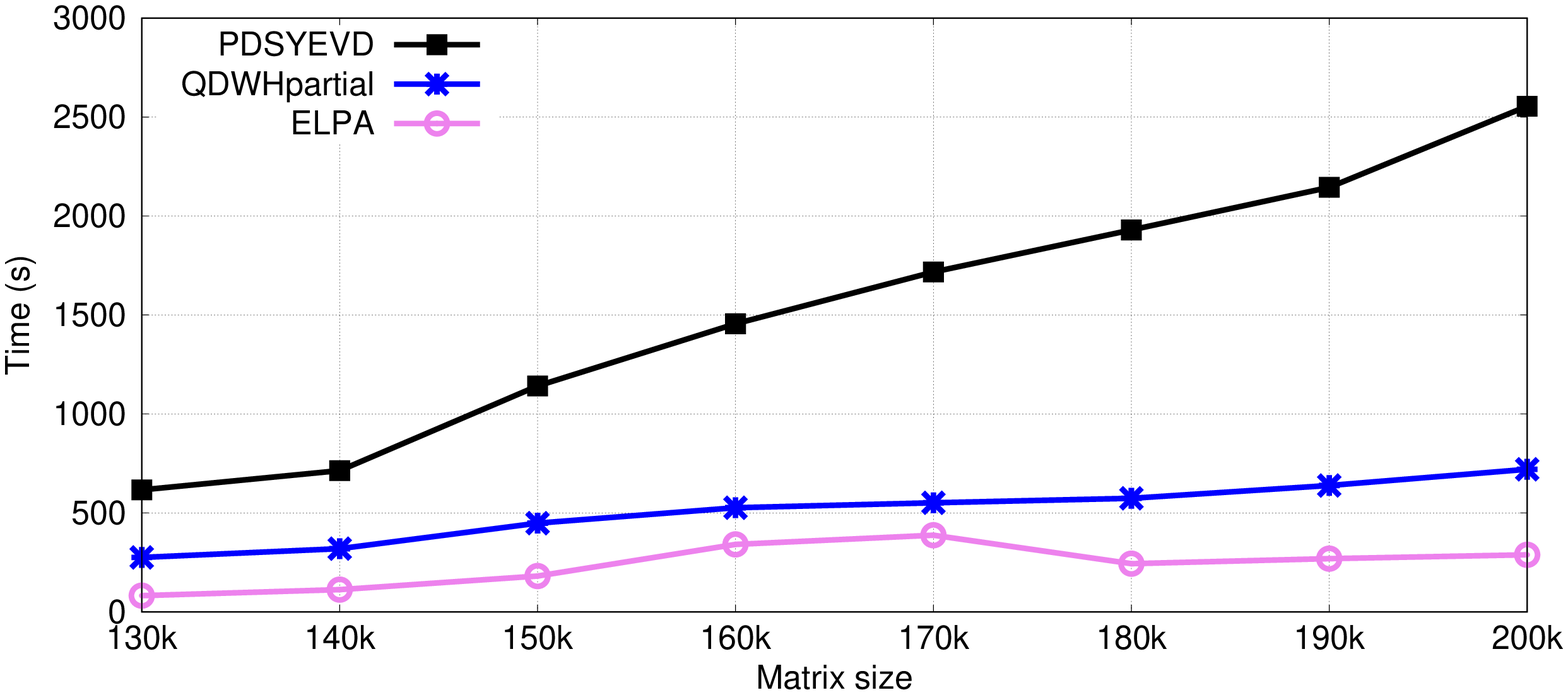}}
      \end{minipage}
      \caption{Performance comparisons in seconds of EIG solvers for various grid topologies.}
      \label{fig:evd-comparison-time}
\end{figure}

Fig.~\ref{fig:evd-comparison-time} reports the performance comparisons of
\texttt{QDWHpartial-EIG} against the two other EIG solvers,
\texttt{ScaLAPACK} \emph{PDSYEVD} and \texttt{ELPA} across various matrix 
sizes and process grid configurations.  \texttt{QDWHpartial-EIG} 
achieves performance superiority as we increase the matrix sizes
up to $1.5X$ on $16 \times 36$, $1.3X$ on $32 \times 72$,
$1.6X$ on $64 \times 144$ and $3.5X$ on $128 \times 288$ grid topologies against
\texttt{ScaLAPACK} \emph{PDSYEVD}. \texttt{QDWHpartial-EIG} remains slower than 
the two-stage approach of \texttt{ELPA}. However, \texttt{QDWHpartial-EIG} 
exposes more parallelism throughout the execution than \texttt{ELPA} 
(i.e., the reduction from band to tridiagonal form is limited in parallelism).
As we increase the number of processors, \texttt{QDWHpartial-EIG} is more capable of
extracting performance from the underlying hardware architecture than 
\texttt{ELPA}.

\begin{figure}[htpb]
      \centering
      \begin{minipage}[htpb]{.46\textwidth}
          \subfigure[P=16 and Q=36.]{
                         \includegraphics[width=.92\linewidth]{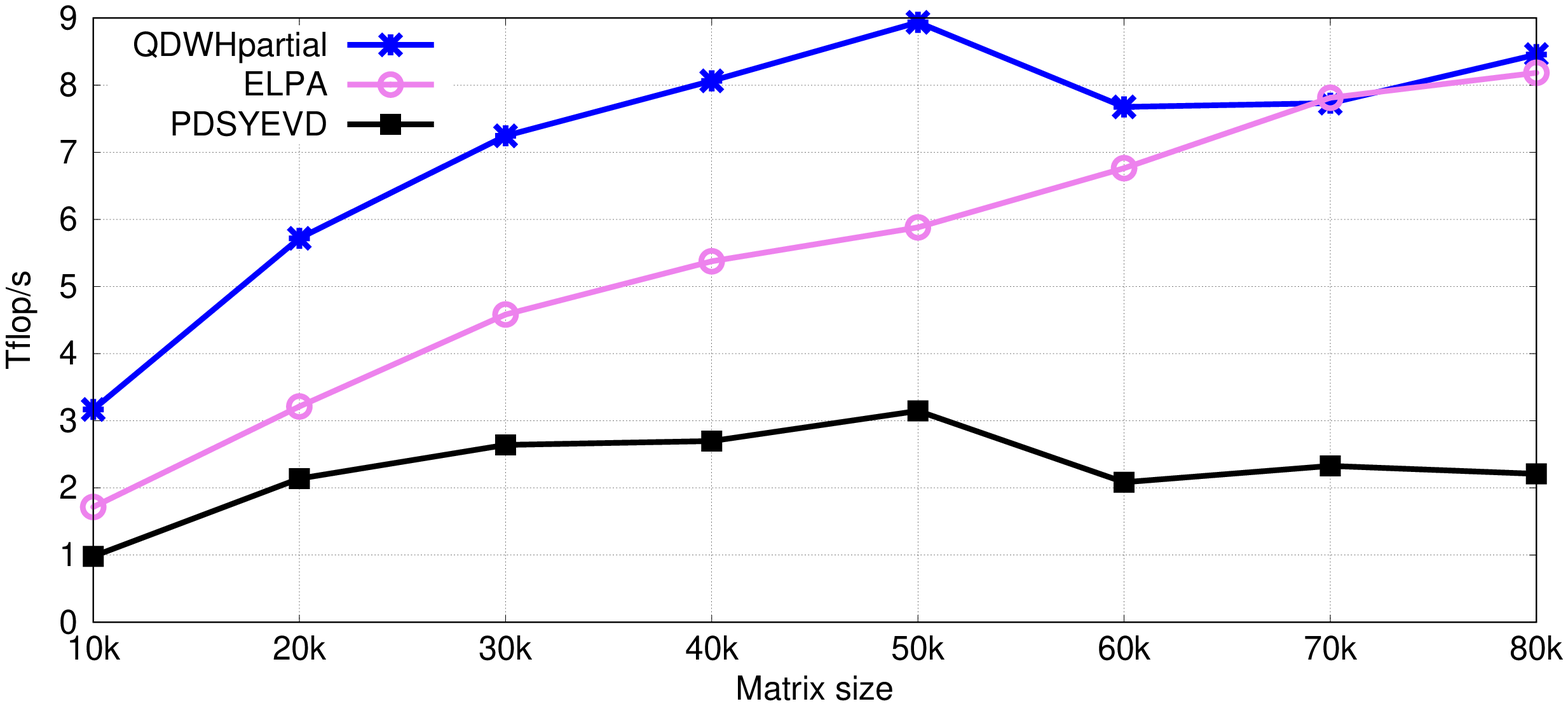}}
          \subfigure[P=32 and Q=72.]{
                         \includegraphics[width=.92\linewidth]{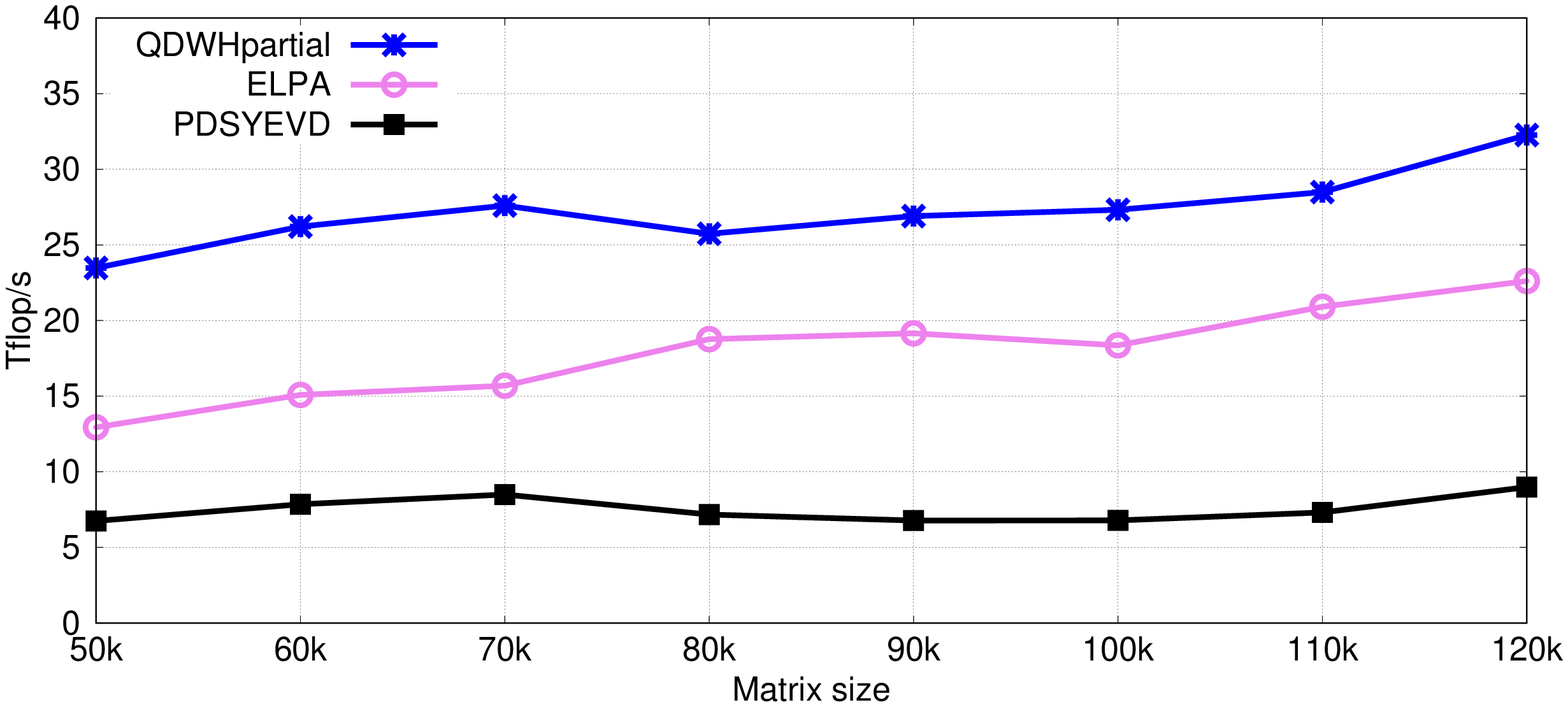}}
      \end{minipage}
      \begin{minipage}[htpb]{.46\textwidth}
          \subfigure[P=64 and Q=144.]{
                         \includegraphics[width=.92\linewidth]{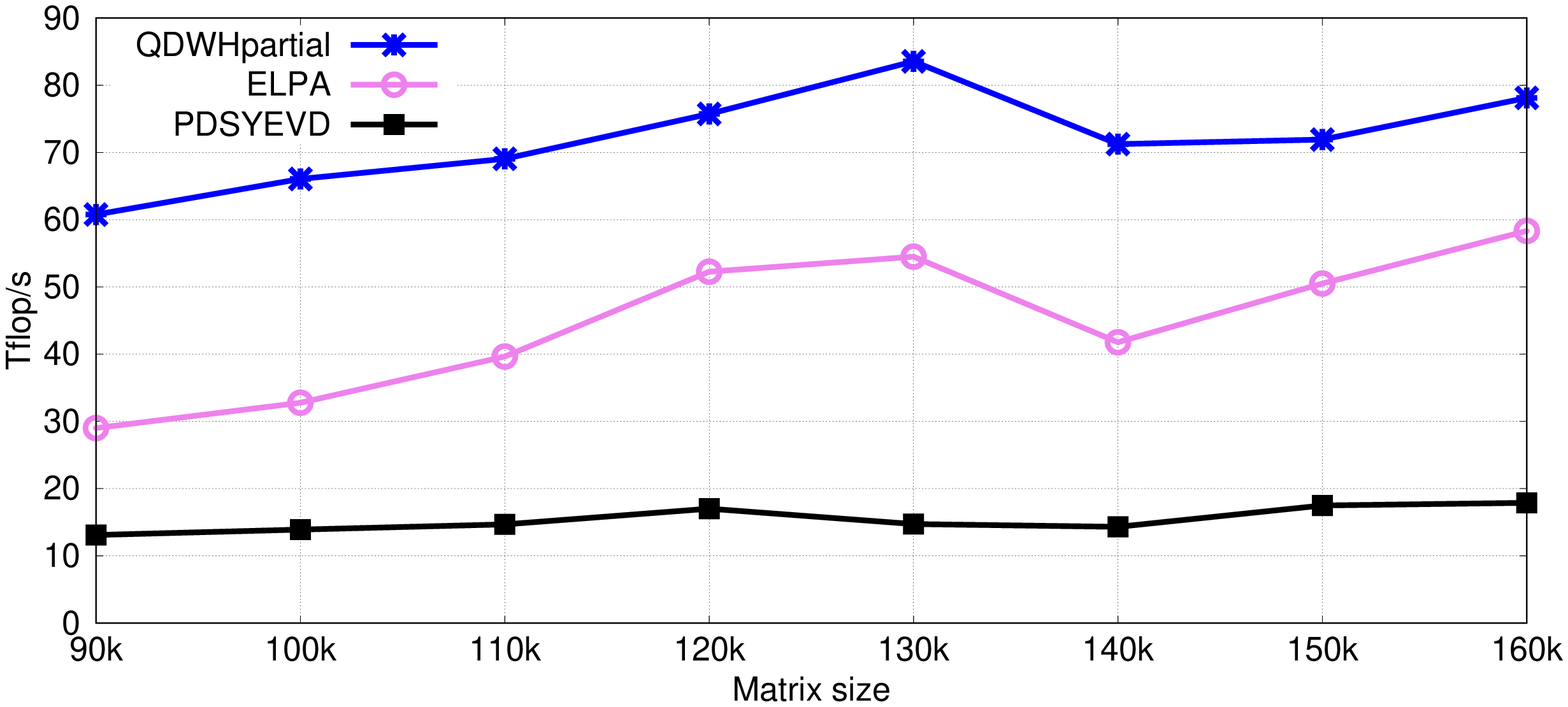}}
          \subfigure[P=128 and Q=288.]{
                         \includegraphics[width=.92\linewidth]{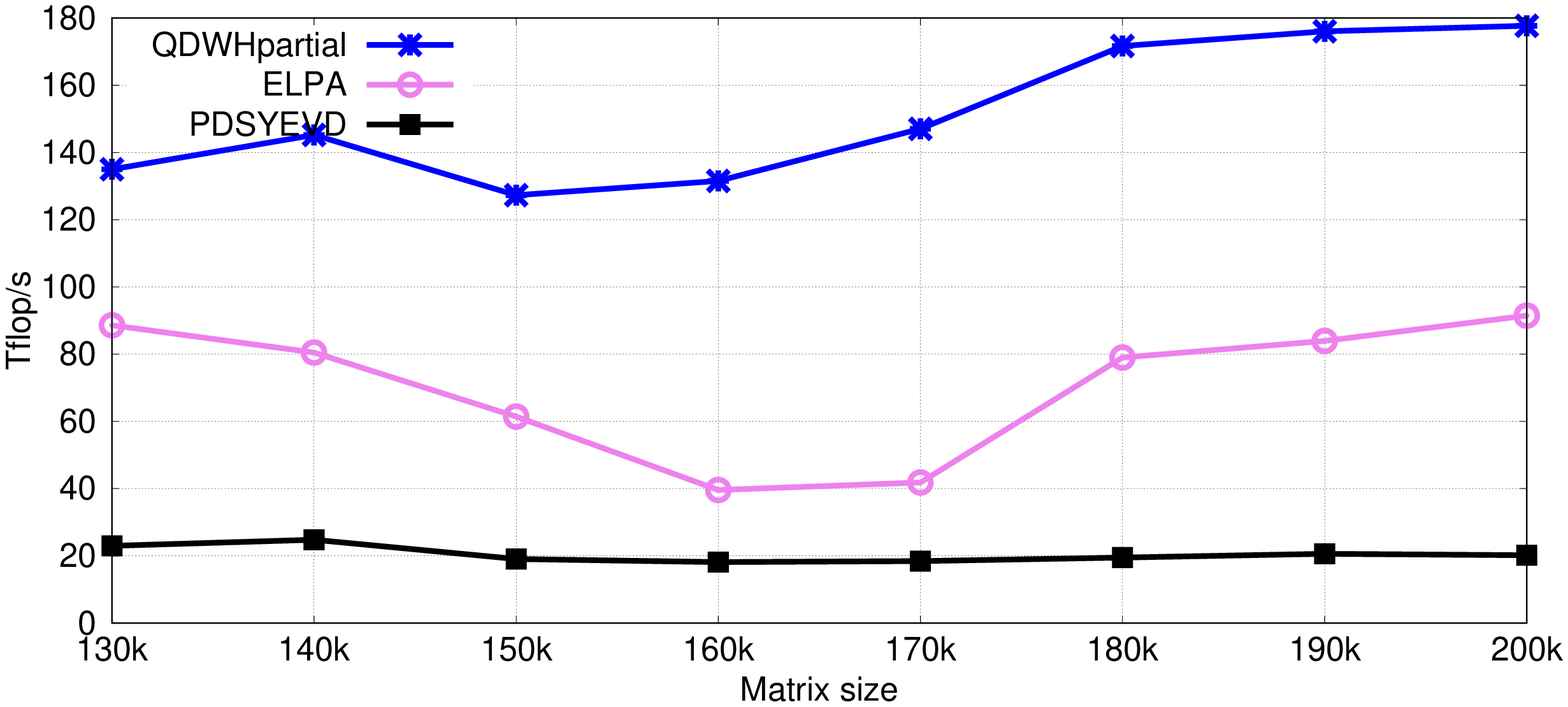}}
      \end{minipage}
      \caption{Performance comparisons in Tflops/s of symmetric EIG solvers for various grid topologies.}
      \label{fig:evd-comparison-flops}
\end{figure}

Indeed, Fig.~\ref{fig:evd-comparison-flops} shows the sustained performance in Tflops/s
and explains why the the performance gap 
between \texttt{QDWHpartial-EIG} and \texttt{ELPA} gets narrower. As we increase the
matrix sizes and the process grids, \texttt{QDWHpartial-EIG} obtains up to a twice higher
rate of executions than \texttt{ELPA}. This performance efficiency may become an advantage 
moving forward with a hardware landscape oriented toward massively parallel
resources delivering high rate of executions (e.g., accelerator-based supercomputers).  

\begin{figure}[htpb]
      \centering
      \begin{minipage}[htpb]{.46\linewidth}
          \subfigure[\texttt{QDWHpartial-SVD} with 13\%.]{
             \includegraphics[width=0.92\linewidth]{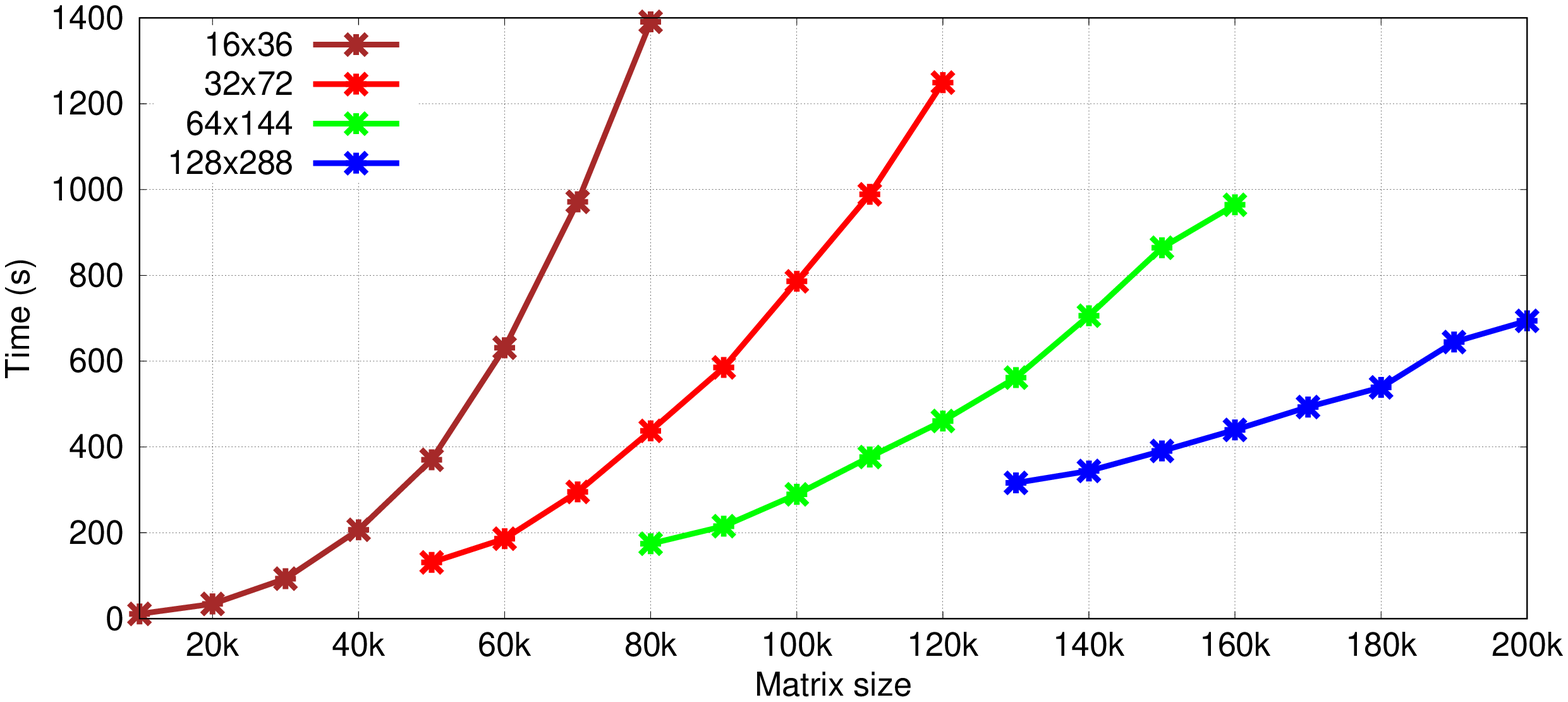}}
      \end{minipage}
      \begin{minipage}[htpb]{.46\linewidth}
          \subfigure[\texttt{QDWHpartial-EIG} with 10\%.]{
             \includegraphics[width=0.92\linewidth]{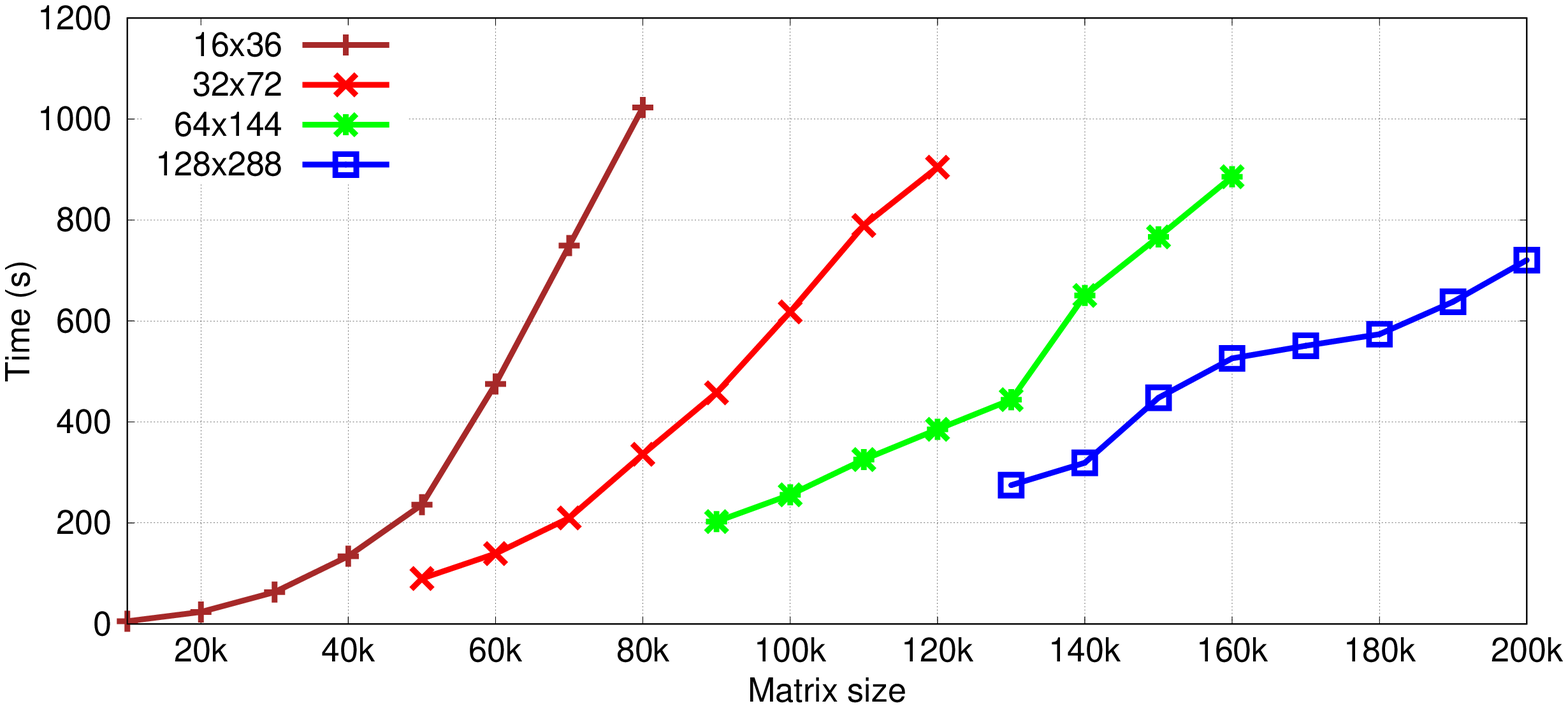}}
      \end{minipage}
	  \caption{Performance scalability of \texttt{QDWHpartial-SVD} and \texttt{QDWHpartial-EIG} for various grid topologies.}
\label{fig:evd-svd-scalability}
\end{figure}

Figure~\ref{fig:evd-svd-scalability} shows various grid topologies
and indicates a decent performance scalability of \texttt{QDWHpartial-SVD} and 
\texttt{QDWHpartial-EIG}, as the matrix sizes increases.
Notice also the various slopes flatten for both solvers, 
since the critical computational phase, i.e., the QDWH-based polar decomposition,
enters into the compute-bound regime of operations along with a better hardware occupancy.
It is also noteworthy to emphasize that the size of the reduced problem may
sometimes be higher than the number of eigenvalues or singular values requested.
This situation explains why \texttt{QDWHpartial-EIG} is sometimes faster than 
\texttt{QDWHpartial-SVD}, although their algorithmic complexities are comparable
(see Table~\ref{table:complexity} in Section~\ref{sec:complexity}).

\section{Summary and future work}
\label{sec:conclusion}
This paper introduces a new algorithm for computing a partial spectrum for the
dense symmetric EIG and SVD solvers. By relying on QDWH-based polar decomposition,
we demonstrate the numerical robustness of \texttt{QDWHpartial-SVD} and \texttt{QDWHpartial-EIG}
against their counterpart routines from state-of-the-art open-source (i.e., \texttt{ELPA} and \texttt{KSVD}) 
and vendor-optimized (i.e., \texttt{ScaLAPACK} from Cray Scientific Library) numerical libraries. 
While \texttt{QDWHpartial-SVD} outperforms the existing approaches up to $6X$, \texttt{QDWHpartial-EIG} shows
performance superiority up to $3.5X$ against the one-stage approach of \emph{PDSYEVD} from \texttt{ScaLAPACK}
but remains slower compared to \texttt{ELPA}. We believe that the inherent massively parallel and 
compute-bound regime of \texttt{QDWHpartial-EIG} may help in narrowing the performance gap observed
against \texttt{ELPA} moving forward with hardware rich in concurrency. We plan to further improve 
our current implementation by using ZOLO-based polar decomposition~\cite{ltaief2019}. Recent
work on leveraging task-based programming model associated with dynamic runtime systems for 
tackling heterogeneous hardware environment~\cite{Sukkaritpds} may also be considered to further speed up the
current implementation on distributed-memory systems equipped with GPU accelerators.


\section*{Acknowledgment}
The authors would like to thank 
Cray Inc. and Intel in the context of the Cray Center
of Excellence and Intel Parallel Computing Center
awarded to ECRC at KAUST.
For computer time, this research used \emph{Shaheen-2} supercomputer hosted 
at the Supercomputing Laboratory at KAUST.

\bibliographystyle{unsrt}
\bibliography{qdwhpartial}
\end{document}